\documentclass{article}
\usepackage[utf8]{inputenc}
\usepackage{amsfonts, amsmath, amssymb, amsthm, authblk, graphicx, enumitem}
\usepackage{cite}
\usepackage[colorlinks=true,urlcolor=blue, citecolor=red,linkcolor=blue]{hyperref}
\usepackage[left=2cm, right=2cm, top=3cm]{geometry}
\usepackage[capitalise]{cleveref}
\newtheorem{theorem}{Theorem}

\newtheorem{definition}[theorem]{Definition}

\newtheorem{lemma}[theorem]{Lemma}
\newtheorem{proposition}[theorem]{Proposition}
\newtheorem{remark}[theorem]{Remark}
\usepackage{dsfont}
\usepackage{hyperref}

\renewcommand{\d}{\displaystyle}
\usepackage{graphicx}
\usepackage{apptools}
\AtAppendix{\counterwithin{lem}{section}}
\newtheorem{lem}{Lemma}
\AtAppendix{\counterwithin{thm}{section}}
\newtheorem{thm}{Theorem}
\newcommand{\pts}[1]{\left(#1\right)}  
\newcommand{\cts}[1]{\left[#1\right]}                                  	  
\newcommand{\lvs}[1]{\left\{#1\right\}}                                	  
\newcommand{\abs}[1]{\left|#1\right|}                                  	  

\newcommand{\C}{\mathbb{C}} 
\newcommand{\R}{\mathbb{R}}
\newcommand{\A}{\mathcal{A}}
\newcommand{\al}{\alpha}

\newcommand{\U}{\mathcal{U}}
\newcommand{\OO}{\mathcal{O}}
\newcommand{\Fi}{\Phi}
\newcommand{\D}{C_c^\infty(0,1)}
\newcommand{\HH}{\mathcal{H}}

\begin{document}
\title{Boundary controllability for a 1D degenerate parabolic equation with drift and a singular potential and a Neumann boundary condition
}

\author[1]{Leandro Galo-Mendoza\thanks{jesus.galo@im.unam.mx}}
\affil[1,2]{Unidad Cuernavaca, Instituto de Matem\'aticas, Universidad Nacional Aut\'onoma de
M\'exico, M\'exico.}
\author[2]{Marcos L\'opez-Garc\'ia\thanks{marcos.lopez@im.unam.mx}}

\maketitle

\abstract{We prove the null controllability of a one-dimensional degenerate parabolic equation with drift and a singular potential. Here, we consider a weighted Neumann boundary control at the left endpoint, where the potential arises. We use a spectral decomposition of a suitable operator, defined in a weighted Sobolev space, and the moment method by Fattorini and Russell to obtain an upper estimate of the cost of controllability. We also obtain a lower estimate of the cost of controllability by using a representation theorem for analytic functions of exponential type.}

\section{Introduction and main results}
Let $T > 0$ and set $Q := (0, 1) \times (0, T )$. For $\al,\beta\in \R$ with $0\leq\al<2$, $\al+\beta >1$, consider the system
\begin{equation}\label{problem}
\left\{\begin{aligned}
u_t-(x^\al u_x)_x-\beta x^{\al -1}u_x-\frac{\mu}{x^{2-\al}} u&=0 & & \text { in }Q, \\
\pts{x^{-\gamma}u_x}(0, t) =f(t),\quad u(1, t)&=0 & & \text { on }(0, T), \\
u(x, 0) &=u_{0}(x) & & \text { in }(0, 1),
\end{aligned}\right.
\end{equation}
provided that $\mu\in\R$ satisfies
\begin{equation}\label{mucon}
-\infty<\mu<\mu(\al+\beta),
\end{equation}
where
\begin{equation}\label{gamadef}
\mu(\delta):=\frac{(\delta-1)^2}{4}, \quad \delta\in\R,\quad \text{and}\quad\gamma=\gamma(\alpha, \beta,\mu):=-(1+\alpha+\beta)/2-\sqrt{\mu(\alpha+\beta)-\mu}.
\end{equation}

The first goal of this work is to provide a notion of a weak solution for the system (\ref{problem}) and show the well-posedness of this problem in suitable interpolation spaces. Here we consider a weighted Neumann boundary condition at the left endpoint to compensate for the singularity of the potential at this point. Then, we use the moment method introduced by Fattorini and Russell in \cite{Fatorini} to prove the null controllability and show an upper bound estimate of the cost of controllability. Next, we use a representation theorem for analytic functions of exponential type to get a lower bound estimate of the cost of controllability.\\

In particular, when $\beta=0$ this work solves the case of strong degeneracy with singularity. Concerning the strongly degenerate case ($1<\al<2$) with no singularity ($\mu=0$), in \cite{gueye2} the authors study the null controllability of a degenerate parabolic equation with a degenerate one-order transport term. In \cite{du,flores,flores2} the authors prove the null controllability of 1D degenerate parabolic equations with first-order terms by means of Carleman inequalities, so they use interior controls.\\ 

Now, assume the system (\ref{problem}) admits a unique solution for initial conditions in a certain Hilbert space $H$, which is described in the next section. We say that the system (\ref{problem}) in null controllable in $H$ at time $T>0$ with controls in $L^2(0,T)$ if for any $u_0\in H$ there exists $f\in L^2(0,T)$ such that the corresponding solution satisfies $u(\cdot,T)\equiv 0$.\\

Once we know the system (\ref{problem}) is null controllable we study the behavior of the cost of the controllability. Consider the set of admissible controls 
\[U(T,\al,\beta,\mu,u_0)=\{f\in L^2(0,T): u \text{ is solution of the system (\ref{problem}) that satisfies } u(\cdot,T)\equiv 0\}.\]
Then the cost of the controllability is defined as
\[\mathcal{K}(T,\al,\beta,\mu):=\sup_{\|u_0\|_H\leq 1}\inf\{\|f\|_{L^2(0,T)}:f\in U(T,\al,\beta,\mu,u_0)\}.\]

In \cite{GaloLopez} it was proved the null controllability of the system (\ref{problem}) provided that $0\leq \al<2$, $\al+\beta<1$, $\mu<\mu(\al+\beta)$, and considering suitable weighted Dirichlet boundary condition at the left endpoint. The main result of this work considers the case $\al+\beta>1$:
\begin{theorem}\label{Teo1}
	Let $T>0$ and $\al,\beta,\mu,\gamma\in \R$ with $0\leq\al<2$, $\al+\beta >1$, $\mu$ and $\gamma$ satisfying (\ref{mucon}) and (\ref{gamadef}) respectively. The next statements hold.
	\begin{enumerate}
		\item \textbf{Existence of a control} For any $u_0\in L^2((0,1);x^\beta dx)$ there exists a control $f \in L^2(0, T )$ such that the solution $u$ to (\ref{problem}) satisfies $u(\cdot,T ) \equiv 0$.
		\item \textbf{Upper bound of the cost} There exists a constant $c>0$ such that for every $\delta\in (0,1)$ we have
		\[\mathcal{K}(T,\alpha,\beta,\mu)\leq c M(T,\alpha,\nu,\delta)T^{1/2}\kappa_\al^{-1/2}
		\exp\pts{-\frac{T}{2}\kappa_\alpha^2 j_{\nu,1}^2},\]
		where 
		\begin{equation}\label{Nu}
			\kappa_\al:=\frac{2-\al}{2},\quad \nu=\nu(\al,\beta,\mu):=\sqrt{\mu(\al+\beta)-\mu}/\kappa_\al,
		\end{equation}
		$j_{\nu,1}$ is the first positive zero of the Bessel function $J_\nu$ (defined in the Appendix), and
		\[M(T,\alpha,\nu,\delta)=\pts{1+\frac{1}{(1-\delta)\kappa_\alpha^2 T}}\cts{\exp\pts{\frac{1}{\sqrt{2}\kappa_\alpha}}+\frac{1}{\delta^3}\exp\pts{\frac{3}{(1-\delta)\kappa_\alpha^2 T}}}\exp\pts{-\frac{(1-\delta)^{3/2}T^{3/2}}{8(1+T)^{1/2}}\kappa_\alpha^3 j_{\nu,1}^2}.\]
		\item \textbf{Lower bound of the cost} There exists a constant $c>0$ such that
		\[\frac{c2^{\nu} \Gamma(\nu+1) \left|J_{\nu}^{\prime}\left(j_{\nu, 1}\right)\right|\exp{\left(\pts{\frac{1}{2}-\frac{\log 2}{\pi}}j_{\nu,2}\right)}}{\pts{{2T \kappa_\alpha}}^{1/2}\left(j_{\nu, 1}\right)^{\nu}}\exp\pts{-\pts{j_{\nu,1}^2+\frac{j_{\nu,2}^2}{2}}\kappa_\alpha^2 T}\leq \mathcal{K}(T,\alpha,\beta,\mu),\]
	\end{enumerate}
	where $j_{\nu,2}$ is the second positive zero of the Bessel function $J_\nu$.
	\end{theorem}

To prove this result we proceed as in \cite{GaloLopez}, in particular, we use the biorthogonal family $(\psi_k)_k$ defined in (\ref{psieta}) and constructed in \cite{GaloLopez}. We also exploit this approach to show the null controllability of the system when the control is located at the right endpoint. Hence, consider the following system
\begin{equation}\label{problem2}
	\left\{\begin{aligned}
		u_t-(x^\al u_x)_x-\beta x^{\al -1}u_x-\frac{\mu}{x^{2-\al}} u&=0 & & \text { in }Q, \\
		\pts{x^{-\gamma}u_x}(0, t) = 0, u(1, t)& = f(t) & & \text { on }(0, T), \\
		u(x, 0) &=u_{0}(x) & & \text { in }(0, 1),
	\end{aligned}\right.
\end{equation}
the corresponding set of admissible controls
\[\widetilde{U}(T,\al,\beta,\mu,u_0)=\{f\in L^2(0,T): u \text{ is solution of the system (\ref{problem2}) that satisfies } u(\cdot,T)\equiv 0\}.\]
and the cost of the controllability given by
\[\widetilde{\mathcal{K}}(T,\al,\beta,\mu):=\sup_{\|u_0\|_H\leq 1}\inf\{\|f\|_{L^2(0,T)}:f\in \widetilde{U}(T,\al,\beta,\mu,u_0)\}.\]
\begin{theorem}\label{Teo2}
	Let $T>0$ and $\al,\beta,\mu,\gamma\in \R$ with $0\leq\al<2$, $\al+\beta >1$, $\mu$ and $\gamma$ satisfying (\ref{mucon}) and (\ref{gamadef}) respectively. The next statements hold.
	\begin{enumerate}
		\item \textbf{Existence of a control} For any $u_0\in L^2((0,1);x^\beta dx)$ there exists a control $f \in L^2(0, T )$ such that the solution $u$ to (\ref{problem2}) satisfies $u(\cdot,T ) \equiv 0$.
		\item \textbf{Upper bound of the cost} There exists a constant $c>0$ such that for every $\delta\in (0,1)$ we have
		\[\widetilde{\mathcal{K}}(T,\alpha,\beta,\mu)\leq \frac{c M(T,\alpha,\nu,\delta)T^{1/2}}{\pts{2\kappa_\al}^{\nu}\Gamma(\nu+1)}
		\pts{\dfrac{2\nu+1}{T}}^{(2\nu+1)/4} \exp\pts{-\frac{2\nu+1}{4}}\exp\pts{-\frac{T}{4}\kappa_\alpha^2 j_{\nu,1}^2}.\]
		\item \textbf{Lower bound of the cost} There exists a constant $c>0$ such that
		\[\frac{c\exp{\left(\pts{\frac{1}{2}-\frac{\log 2}{\pi}}j_{\nu,2}\right)}}{T^{1/2} \kappa_\alpha^{3/2}j_{\nu, 1}}\exp\pts{-\pts{j_{\nu,1}^2+\frac{j_{\nu,2}^2}{2}}\kappa_\alpha^2 T}\leq \widetilde{\mathcal{K}}(T,\alpha,\beta,\mu).\]
	\end{enumerate}
	\end{theorem}

Finally, we also analyze the null controllability of the system when the parameters satisfy $0\leq \al<2$, $\beta=1-\al$, and $\mu<0$. Thus, we consider the following system.
\begin{equation}\label{problem3}
	\left\{\begin{aligned}
		u_t-(x^\al u_x)_x-(1-\al) x^{\al -1}u_x-\frac{\mu}{x^{2-\al}} u&=0 & & \text { in }Q, \\
		\pts{x^{\sqrt{-\mu}}u}(0, t) = f(t), u(1, t)& = 0 & & \text { on }(0, T), \\
		u(x, 0) &=u_{0}(x) & & \text { in }(0, 1),
	\end{aligned}\right.
\end{equation}
The corresponding set of admissible controls is given by
\[\widehat{U}(T,\al,\mu,u_0)=\{f\in L^2(0,T): u \text{ is solution of the system (\ref{problem3}) that satisfies } u(\cdot,T)\equiv 0\},\]
and the cost of the controllability is given by
\[\widehat{\mathcal{K}}(T,\al,\mu):=\sup_{\|u_0\|_H\leq 1}\inf\{\|f\|_{L^2(0,T)}:f\in \widehat{U}(T,\al,\mu,u_0)\}.\]
We use some result from the singular Sturm-Liouville theory to show the well-posedness of system (\ref{problem3}).
\begin{theorem}\label{unosuma}
Let $T>0$ and $\al,\mu\in \R$ with $0\leq\al<2$, $\mu<0$. The next statements hold.
\begin{enumerate}
\item \textbf{Existence of a control} For any $f \in L^2(0,T)$ and $u_0\in L^2((0,1);x^{1-\al} dx)$ there exists a control $f \in L^2(0, T)$ such that the solution $u$ to (\ref{problem3}) satisfies $u(\cdot,T ) \equiv 0$.
\item \textbf{Upper bound of the cost} There exists a constant $c>0$ such that for every $\delta\in (0,1)$ we have
		\[\widehat{\mathcal{K}}(T,\alpha,\mu)\leq \frac{c M(T,\alpha,\nu,\delta)T^{1/2}}{\kappa_\al^{1/2}\sqrt{-\mu}}
\exp\pts{-\frac{T}{2}\kappa_\alpha^2 j_{\nu,1}^2},\]
where $\nu=\nu(\al,\mu):=\sqrt{-\mu}/\kappa_\al$.
\item \textbf{Lower bound of the cost} There exists a constant $c>0$ such that
		\[\frac{c2^{\nu} \Gamma(\nu+1) \left|J_{\nu}^{\prime}\left(j_{\nu, 1}\right)\right|\exp{\left(\pts{\frac{1}{2}-\frac{\log 2}{\pi}}j_{\nu,2}\right)}}{\pts{{2T \kappa_\alpha}}^{1/2}\sqrt{-\mu}\left(j_{\nu, 1}\right)^{\nu}}\exp\pts{-\pts{j_{\nu,1}^2+\frac{j_{\nu,2}^2}{2}}\kappa_\alpha^2 T}\leq \widehat{\mathcal{K}}(T,\alpha,\mu).\]
\end{enumerate}
\end{theorem}

This paper is organized as follows. In Section \ref{setting}, we introduce suitable weighted Sobolev spaces and prove some results about the trace (at the endpoints) of functions in these spaces, as well as on the behavior of these functions at the endpoints, we also show an integration by parts formula. In that section, we prove that the autonomous operator given in (\ref{Aoperator}) is diagonalizable, which allows the introduction of interpolation spaces for the initial data. Then, we prove the system (\ref{problem}) is well-posed in this setting. \\

In Section \ref{leftend} we prove Theorem \ref{Teo1} by using the moment method, as a consequence, we get an upper estimate of $\mathcal{K}(T,\al,\beta,\mu)$. Then we use the representation theorem in Theorem \ref{repreTheo} to obtain a lower estimate of $\mathcal{K}(T,\al,\beta,\mu)$. In Section \ref{rightend} we proceed as before to prove  Theorem \ref{Teo2}. Finally, in Section \ref{sumauno} we sketch the proof of Theorem 3.

\section{Functional setting and well-posedness}\label{setting}
In this section, we introduce some suitable weighted spaces. First, consider the weighted Lebesgue space $L^2_\beta(0,1):=L^2((0,1);x^\beta dx)$, $\beta\in\R$, endowed with the inner product
\[\langle f,g\rangle_\beta:=\int_0^1 f(x)g(x)x^\beta dx,\]
and its corresponding norm denoted by $\|\cdot\|_{\beta}$.\\

For $\al,\beta\in \R$ consider the weighted Sobolev space
\[H_{\alpha, \beta}^{1}(0,1)=\left\{u \in L_{\beta}^{2}(0,1)\cap H^1_{loc}(0,1): x^{\alpha / 2} u_{x} \in L_{\beta}^{2}(0,1)\right\}\]
endowed with the inner product
\[\langle u, v\rangle_{\alpha, \beta}:=\int_{0}^{1} u v\, x^{\beta} \!dx+\int_{0}^{1} x^{\alpha+\beta} u_{x} v_{x} dx,\]
and its corresponding norm denoted by $\|\cdot\|_{\alpha,\beta}$.\\

The next result implies that we can talk about the trace at $x=1$ of functions in $H_{\alpha, \beta}^{1}(0,1)$.
\begin{proposition} \label{abscont}
	Let $\al,\beta\in \R$. Then $H_{\alpha, \beta}^{1} \subset W^{1,1}(\varepsilon,1)$ for all $\varepsilon \in \pts{0,1}$. In particular, $H_{\al,\beta}^{1}(0,1)\subset C((0,1])$, and $|u|^2\in W^{1,1}(\varepsilon,1)$ for all $u\in H_{\alpha, \beta}^{1}(0,1),\varepsilon \in \pts{0,1}$.
\end{proposition}
\begin{proof}
	Let $u\in H_{\alpha, \beta}^{1}(0,1)$. For $\varepsilon\in \pts{0,1},\delta\in\R$ fixed, there exists a constant $c(\varepsilon,\delta)>0$ such that $x^\delta\leq c(\varepsilon,\delta)$, $x\in(\varepsilon,1]$, thus
	\[\int_{\varepsilon}^{1}|u| \mathrm{d}x \leq (1-\varepsilon)^{1/2}\left(\int_{\varepsilon}^{1} |u|^{2}  \mathrm{d}x\right)^{1 / 2}\leq (1-\varepsilon)^{1/2}c(\varepsilon,-\beta)^{1/2}\left(\int_{0}^{1} |u|^{2}  x^{\beta}\mathrm{d}x\right)^{1 / 2},\quad\text{and  }\]
	\[\int_{\varepsilon}^{1}\left|u_{x}\right| \mathrm{d}x\leq (1-\varepsilon)^{1/2}c(\varepsilon,-\al-\beta)^{1/2}\left(\int_{0}^{1} |u_{x}|^{2}  x^{\al+\beta}\mathrm{d}x\right)^{1 / 2}.\]
Hence we get the existence of the limit  $u(1):=\lim_{x\rightarrow 1^-}u(x)$, and $u\in C\pts{\cts{\varepsilon,1}}$.
\end{proof}

\begin{definition}
For $\al,\beta\in \R$ consider the space
\[H_{\alpha, \beta,N}^{1}=H_{\alpha, \beta,N}^{1}(0,1):=\left\{ u \in H_{\alpha, \beta}^{1}(0,1): u(1)=0\right\}.\]
\end{definition}

Next, we generalize the so-called Hardy inequality in the setting of the weighted Sobolev space $H_{\alpha, \beta,N}^{1}$.
\begin{proposition}\label{Hardy}
	For $\al,\beta\in \R$ with $\al+\beta>1$, the Hardy inequality
	\begin{equation}\label{hardy}
		\mu(\al+\beta)\int_0^1\frac{|u|^2}{x^{2-(\al+\beta)}}\mathrm{d}x\leq\int_0^1x^{\al+\beta} |u_x|^2\mathrm{d}x
	\end{equation}
	holds for any $u\in H_{\alpha, \beta,N}^{1}$. In particular, $H_{\alpha, \beta,N}^{1}\hookrightarrow L^2_{\al+\beta-2}(0,1)$.
\end{proposition}
\begin{proof}
	Let $u\in H_{\alpha, \beta,N}^{1}$ and $\varepsilon\in (0,1)$. Set $\delta=\al+\beta$. Since $|u|^2\in W^{1,1}(\varepsilon,1)$ we have
	\begin{eqnarray*}
		\int_\varepsilon^1\pts{x^{\delta/2}u_x-\frac{1-\delta}{2}\frac{u}{x^{(2-\delta)/2}}}^2\mathrm{d}x
		& = & \int_\varepsilon^1x^{\delta}|u_x|^2\mathrm{d}x+\mu(\delta)\int_\varepsilon^1\frac{u^2}{x^{2-\delta}}\mathrm{d}x-\frac{1-\delta}{2}\int_\varepsilon^1\frac{(u^2)_x}{x^{1-\delta}}\mathrm{d}x\\
		& = & \int_\varepsilon^1x^{\delta}|u_x|^2\mathrm{d}x-\mu(\delta)\int_\varepsilon^1\frac{u^2}{x^{2-\delta}}\mathrm{d}x-\frac{1-\delta}{2}\pts{\lim_{x\rightarrow 1^{-}}\frac{|u(x)|^2}{x^{1-\delta}}-\frac{|u(\varepsilon)|^2}{\varepsilon^{1-\delta}}}\\
		& = &
		\int_\varepsilon^1x^{\delta}|u_x|^2\mathrm{d}x-\mu(\delta)\int_\varepsilon^1\frac{u^2}{x^{2-\delta}}\mathrm{d}x+\frac{1-\delta}{2}\frac{|u(\varepsilon)|^2}{\varepsilon^{1-\delta}},		
	\end{eqnarray*}
	since $\delta > 1$ we get
	$$
	\mu(\delta)\int_\varepsilon^1\frac{u^2}{x^{2-\delta}}\mathrm{d}x\leq\int_\varepsilon^1x^\delta u_x^2\mathrm{d}x
	$$
	for all $\varepsilon \in (0,1)$. The result follows by the dominated convergence theorem.  
\end{proof}

The next result will allow us analyze the behavior at $x=0$ of functions in $H_{\al, \beta, N}^{1}$,  see (\ref{vanish0}).
\begin{proposition}\label{gamma_u}
	Let $\al,\beta\in \R$ with $\al+\beta >1$. Then $x^{\delta}u\in W^{1,1}(0,1)$ for all $u\in H_{\al, \beta, N}^{1}$ provided that $\delta > (\al+\beta-1)/2$. 
\end{proposition}
\begin{proof} Let $u\in H_{\al, \beta, N}^{1}$ and assume $2\delta > \al+\beta-1$. 

	Since $\pts{x^{\delta}u}_{x} = x^{\delta}u_{x}+\delta x^{\delta-1}u$, we compute 
	\begin{equation}\label{apoyo}
	\int_{0}^{1}x^{\delta}|u_{x}| \mathrm{d}x \leq \frac{1}{(2\delta-(\al+\beta) +1)^{1/2}}\left(\int_{0}^{1} x^{\al+\beta}|u_{x}|^{2} \mathrm{d}x\right)^{1 / 2}<\infty,
	\end{equation}
	and Proposition \ref{Hardy} implies
	\begin{eqnarray*}
	\int_{0}^{1}x^{\delta-1}|u| \mathrm{d}x &\leq& \frac{1}{(2\delta-(\al+\beta) +1)^{1/2}}\left(\int_{0}^{1}\frac{|u|^{2}}{x^{2-(\al+\beta)}}\mathrm{d}x\right)^{1 / 2}\\
	&\leq&\frac{1}{(2\delta-(\al+\beta) +1)^{1/2}}\frac{1}{\mu(\al+\beta)}\left(\int_{0}^{1}x^{\al+\beta}u_{x}^{2}\mathrm{d}x\right)^{1 / 2}<\infty.
	\end{eqnarray*}
	Hence $\pts{x^{\delta}u}_{x}\in L^1(0,1)$. Notice that $x^\delta\leq x^{\delta -1}$ on $(0,1)$, thus $x^{\delta}u\in L^1(0,1)$ and the result follows.
	\end{proof}
\begin{remark}\label{cerolim}
	For $\al,\beta\in \R$ with $\al+\beta>1$, the last result implies the existence of $L_{\delta}:=\lim_{x\to 0^{+}}x^{\delta}u(x)$ provided that $\delta> (\al + \beta - 1)/2$, in fact, $L_{\delta} = 0$. Now choose any $\delta> (\al + \beta - 1)/2$ so
\begin{eqnarray*}
x^{\delta}|u(x)|&\leq&\int_0^x \left|\frac{\mathrm{d}}{\mathrm{d}s} (s^\delta u(s))\right|\mathrm{d}s\\
			&\leq&\frac{x^{\delta-(\al+\beta-1)/2}}{(2\delta-(\al+\beta)+1)^{1/2}}\cts{\pts{\int_{0}^{x}s^{\al+\beta}|u_x|^2\mathrm{d}s}^{1/2}+\delta\pts{\int_{0}^{x}\frac{|u|^2}{s^{2-(\al+\beta)}}\mathrm{d}s}^{1/2}},
\end{eqnarray*}
therefore
\begin{equation}\label{vanish0}
\lim_{x\to 0^+}x^{(\al+\beta-1)/2}|u(x)|=0, \quad u\in H_{\alpha, \beta,N}^{1}.
\end{equation}
\end{remark}

From now on we assume $\al < 2$ and $\al + \beta>1$. For any $u \in H_{\alpha, \beta, N}^{1}$, we obtain the weighted Poincaré inequality from Proposition \ref{Hardy}:
\begin{equation}\label{poincare}
	\int_{0}^{1}x^{\beta}|u|^{2}\mathrm{d}x \leq\int_{0}^{1}\dfrac{|u|^2}{x^{2-(\al+\beta)}} \mathrm{d}x\leq \frac{1}{\mu(\al+\beta)}\int_{0}^{1} x^{\alpha+\beta}\left|u_{x}\right|^{2} \mathrm{d}x,
\end{equation}
therefore
\[\|u\|_{\al,\beta,N}:=\left(\int_{0}^{1} x^{\alpha+\beta}\left|u_{x}\right|^{2} \mathrm{d}x\right)^{1 / 2}\]
is an equivalent norm to $\|u\|_{\alpha, \beta}$ in $H_{\al, \beta, N}^{1}$.\\

For $\mu<\mu(\alpha+\beta)$, Proposition \ref{Hardy} also implies that
\[
\|u\|_{*}=\left(\int_{0}^{1} x^{\al+\beta}\left[\left|u_{x}\right|^{2}-\frac{\mu}{x^{2}} u^{2}\right] \mathrm{d}x\right)^{1 / 2}
\]
is an equivalent norm to $\|u\|_{\al,\beta,N}$ in $H_{\alpha, \beta, N}^{1}$. We have
\begin{eqnarray*}
	\|u\|_{\al,\beta,N} &\leq & \| u\|_{*} \leq \left(1-\frac{\mu}{\mu(\alpha+\beta)}\right)^{1/2}\| u \|_{\al,\beta,N},\quad \mu<0, \\
	\left(1-\frac{\mu}{\mu(\alpha+\beta)}\right)^{1/2}\|u\|_{\al,\beta,N} &\leq & \| u\|_{*} \leq \left\|u \right\|_{\al,\beta,N},\quad 0\leq \mu <\mu(\al+\beta).
\end{eqnarray*}

Since $\D\subset H_{\alpha, \beta, N}^{1}\subset L^2_\beta(0,1),$ and (\ref{poincare}) implies that the inclusion $(H_{\alpha, \beta,N}^1,\|\cdot\|_*)\hookrightarrow L^2_\beta(0,1)$ is continuous, the following definition makes sense.

\begin{definition}
	For $\al,\beta\in \R$ with $\al<2$ and $\al+\beta >1$, consider the Gelfand triple $\left((H_{\alpha, \beta,N}^1,\|\cdot\|_*), L^2_\beta(0,1), H_{\alpha, \beta,N}^{-1}\right)$, i.e $H_{\alpha, \beta,N}^{-1}$ stands for the dual space of $(H_{\alpha, \beta,N}^1,\|\cdot\|_*)$ with respect to the pivot space $L^2_\beta(0,1)$:
	\[(H_{\alpha, \beta,N}^1,\|\cdot\|_*)\hookrightarrow L^2_\beta(0,1)=\left(L^2_\beta(0,1)\right)'\hookrightarrow H_{\alpha, \beta,N}^{-1}:=\pts{H_{\alpha, \beta,N}^1,\|\cdot\|_*}'.\]
\end{definition}

The inner product $\langle \cdot ,\cdot \rangle_*$ induces an isomorphism $\A:H_{\alpha, \beta,N}^1\rightarrow H_{\alpha, \beta,N}^{-1}$ given by
\[\langle u ,v \rangle_*=\langle \A u,v\rangle_{H_{\alpha, \beta,N}^{-1},H_{\alpha, \beta,N}^{1}},\quad u,v\in H_{\alpha, \beta,N}^{1}.\]

Let $D(\A):=\A^{-1}(L^2_\beta(0,1))=\{u\in H_{\alpha, \beta,N}^{1}: \A u\in L^2_\beta(0,1)\}=\{u\in H_{\alpha, \beta,N}^{1}: \exists f\in L^2_\beta(0,1) \text{ such that }\langle u,v \rangle_*=\langle f,v\rangle_\beta, \text{ for all } v\in H_{\alpha, \beta,N}^{1} \}$.\\

The next result gives a handy characterization of $D(\A)$. It shows the behavior of the derivative of functions in $D(\A)$ at the endpoints, see (\ref{deri1}) and (\ref{vanishderi}), and also provides an integration by parts formula, see (\ref{quasi}).
\begin{proposition} For $\al,\beta,\mu\in \R$ with $0\leq\al<2$, $\al+\beta >1$ and $\mu<\mu(\alpha+\beta)$, we have
	\[D(\A)=\left\{u\in H^1_{\al,\beta,N}\cap H^2_{loc}(0,1): (x^\al u_x)_x+\beta x^{\al -1}u_x+\frac{\mu}{x^{2-\al}} u\in L^2_\beta(0,1)\right\}.\]
\end{proposition}
\begin{proof}
	Let $H$ be the set on the right-hand side, we will show that $D(\A)=H$.\\
	
	Pick $u \in D(\A)$, then there exists $f \in L^{2}_\beta(0,1)$ such that 
	\[\int_{0}^{1} \left(x^{\alpha+\beta} u_{x} v_{x}-\frac{\mu}{x^{2-\alpha-\beta}} u v\right)\mathrm{d}x=\int_0^1 f v x^{\beta} \mathrm{d}x
	\quad \text{for all }v \in H_{\alpha, \beta, N}^{1}.
	\]
	In particular,
	\[
	\int_{0}^{1} x^{\alpha+\beta} u_{x} v_{x}\mathrm{d}x=\int_{0}^{1}\left(f+\frac{\mu}{x^{2-\alpha}} u\right) v x^{\beta} \mathrm{d}x\quad \text {for all  } v \in \D,
	\]
	hence
	\[-\left(x^{\alpha+\beta} u_{x}\right)_{x}=\left(f+\frac{\mu}{x^{2-\alpha}}u\right) x^{\beta} \text { in } \D^{\prime},
	\] 
	which implies
	\[
	\left(x^{\alpha} u_{x}\right)_{x}+\beta x^{\alpha-1} u_{x}+\frac{\mu}{x^{2-\alpha}} u=-f \text { in } \D^{\prime},
	\]
	therefore $u\in H$.\\
	
	Now let  $u \in H$. We claim that $x^{\delta} u_{x} \in W^{1,1}(0,1)$ for all $\delta > \pts{\al+\beta+1}/2$. Just apply (\ref{apoyo}) with $\delta-1$ instead of $\delta$ to get that $x^{\delta-1}u_x\in L^1(0,1)$, in particular $x^{\delta}u_x\in L^1(0,1)$. On the other hand, we have
		\begin{eqnarray*}
		\int_{0}^{1}x^{\delta}|u_{xx}|\mathrm{d}x & \leq & \int_{0}^{1}x^{\delta-(\al+\beta/2)}|(x^{\al}u_{x})_{x}+\beta x^{\al-1}u_{x}+\dfrac{\mu}{x^{2-\al}}u|x^{\beta/2}\mathrm{d}x+(\al+\beta)\int_{0}^{1}x^{\delta-1}|u_{x}|\mathrm{d}x\\
		 & & +\quad|\mu|\int_{0}^{1}x^{\delta-(\al+\beta)/2-1}\dfrac{|u|}{x^{(2-\al-\beta)/2}}\mathrm{d}x\\
		 &\leq & \dfrac{1}{(2\delta-(2\al+\beta)+1)^{1/2}}\pts{\int_{0}^{1}|(x^{\al}u_{x})_{x}+\beta x^{\al-1}u_{x}+\dfrac{\mu}{x^{2-\al}}u|^2 x^{\beta}\mathrm{d}x}^{1/2}\\
		 & & + \frac{\al + \beta}{(2\delta-(\al+\beta+1))^{1/2}}\left(\int_{0}^{1}x^{\al+\beta} |u_{x}|^{2} \mathrm{d}x\right)^{1 / 2} +\frac{|\mu|}{(2\delta-(\al+\beta+1))^{1/2}}\pts{\int_{0}^{1}\frac{ |u|^{2}}{x^{2-(\al+\beta)}} \mathrm{d}x}^{1 / 2}.
		 	\end{eqnarray*}
	   Notice the last quantity is finite by Proposition \ref{Hardy}.\\
    
    Thus, we get the existence of the limit
    \begin{equation}\label{deri1}
     u_{x}(1):=\lim_{x \rightarrow 1^{-}}x^{\delta}u_{x}(x),
     \end{equation}
 and we also have that $\lim_{x\rightarrow 0^+}x^{\delta}u_{x}(x)=0$ provided that $\delta > (\al+\beta+1)/2$, see Remark \ref{cerolim}. As in the proof of (\ref{vanish0}), we can see that 
 \begin{equation}\label{vanishderi}
 \lim_{x\rightarrow 0^+}x^{(\al+\beta+1)/2}u_{x}(x)=0.
 \end{equation}
    
    Now consider any $v\in H^{1}_{\al,\beta,N}$. We claim that $x^{\alpha+\beta} u_{x} v \in W^{1,1}(0,1)$:
    \[\int_0^1x^{\al+\beta}|u_x v|\mathrm{d}x\leq \pts{\int_0^1x^{\al+\beta}|u_x|^2\mathrm{d}x}^{1/2}\pts{\int_0^1x^{\al+\beta}|v|^2\mathrm{d}x}^{1/2}\leq \|u\|_{\al,\beta,N}\|v\|_\beta<\infty,\quad\text{and}\]
    \begin{equation}\label{partes}
    	\left(x^{\alpha+\beta} u_{x} v\right)_{x}= x^{\alpha+\beta} u_{x} v_{x}+x^{\beta}\left(\left(x^{\alpha} u_{x}\right)_x+\beta x^{\alpha-1} u_{x}+\frac{\mu}{x^{2-\alpha}} u\right)v -\frac{\mu}{x^{2-\alpha-\beta}} u v \in L^{1}(0,1).
    \end{equation}
    
    On the other hand, (\ref{vanish0}), (\ref{deri1}) and (\ref{vanishderi}) imply that
    \begin{equation}\label{boundaryC}
    	\lim_{x\rightarrow 0^+}x^{\al+\beta}u_{x}(x)v(x)= 0\quad \text{and}\quad \lim_{x\rightarrow 1^-}x^{\al+\beta}u_{x}(x)v(x) = 0.
    \end{equation}
    
	Thus, from (\ref{partes}) we  get
	\begin{equation}\label{quasi}
		\int_{0}^{1} \left(x^{\alpha+\beta} u_{x} v_{x}-\frac{\mu}{x^{2-\alpha-\beta}} u v\right) \mathrm{d}x=-\int_{0}^{1} x^{\beta}\left(\left(x^{\alpha} u_{x}\right)_{x}+\beta x^{\alpha-1} u_{x}+\frac{\mu}{{x^{2-\alpha}}} u\right) v \mathrm{d}x
	\end{equation}
	for all $u \in H, v \in H_{\alpha, \beta, N}^{1}$. Therefore $u \in D(\A)$.
\end{proof}

For $\al,\beta,\mu\in \R$ with $0\leq\al<2$, $\al+\beta>1$, $\mu<\mu(\al+\beta)$, we consider the unbounded operator $\A:D(\A)\subset L^2_\beta(0,1)\rightarrow L^2_\beta(0,1)$ given by
\begin{equation}\label{Aoperator}
	\A u:=-(x^\al u_x)_x-\beta x^{\al -1}u_x-\frac{\mu}{x^{2-\al}} u.
\end{equation}

From Proposition 9 in \cite[p. 370]{dautray} we have that $\A$ is a closed operator with $D(\A)$ dense in $L^2_\beta(0,1)$. We also have that $\A:(D(\A),\|\cdot\|_{D(\A)})\rightarrow L^2_\beta(0,1)$ is an isomorphism, where
\[\|u\|_{D(\A)}=\|u\|_\beta+\|\A u\|_\beta,\quad u\in D(\A).\]

The next result shows that $\mathcal{A}$ is a diagonalizable self-adjoint operator whose Hilbert basis of eigenfunctions can be written in terms of a Bessel function of the first kind $J_{\nu}$ and its corresponding zeros $j_{\nu,k}$, $k\geq 1$, located in the positive half line. In the appendix, we give some properties of Bessel functions and their zeros.
\begin{proposition}\label{bases}
	$-\A$ is a negative self-adjoint operator. Furthermore, the family
	\begin{equation}\label{Phik}
		\Fi_k(x):=\frac{(2\kappa_\al)^{1/2}}{\abs{J'_\nu\pts{j_{\nu,k}}}}x^{\pts{1-\al-\beta}/2}J_\nu\pts{j_{\nu,k}x^{\kappa_\al}},\quad k\geq 1,
	\end{equation}
	is an orthonormal basis for $L^2_\beta(0,1)$ such that
	\begin{equation}\label{lambdak}
		\mathcal{A}\Fi_k=\lambda_k \Fi_k, \quad \lambda_k=\kappa_\al^2 \pts{j_{\nu,k}}^2,\quad k\geq 1,
	\end{equation}
	where $\nu$ is defined in (\ref{Nu}).
\end{proposition}
\begin{proof}
From (\ref{quasi}) we get that $\A$ is a symmetric operator. Letting $u=v\in D(\A)$ in (\ref{quasi}) and using Proposition \ref{Hardy} we obtain that $-\A\leq 0$.\\

We claim that $\textrm{Ran}(I+\A)=L^2_\beta(0,1)$: Let $f \in L_{\beta}^{2}(0,1)$ be given. Since the inner product $\langle\cdot,\cdot\rangle_\beta+\langle\cdot,\cdot\rangle_*$ is equivalent to $\langle\cdot,\cdot\rangle_{\alpha,\beta}$ in $H_{\alpha, \beta, N}^{1}$ and $f\in (H_{\alpha, \beta, N}^{1},\|\cdot\|_{\alpha,\beta})'$, the Riesz representation theorem implies that there exists a unique $u\in H_{\alpha, \beta, N}^{1}$ such that
\[\int_{0}^{1} u v x^{\beta} \mathrm{d}x+\int_{0}^{1} x^{\alpha+\beta}\left(u_{x} v_{x}-\frac{\mu}{x^{2}} u v\right) \mathrm{d}x=\int_{0}^{1} f v x^{\beta} \mathrm{d}x\]
for all $v\in H_{\alpha, \beta, N}^{1}$. Therefore
$$
u-\left(x^{\alpha} u_{x}\right)_{x}-\beta x^{\alpha-1} u_{x}-\frac{\mu}{x^{2-\alpha}} u=f \text { in } \D^{\prime},
$$
thus $u\in D(\A)$ and $u+\A u=f$.\\

It follows that $-\A$ is $m$-dissipative in $L^2_\beta(0,1)$ and Corollary 2.4.10 in \cite[ p. 24]{cazenave} implies that $-\A$ is self-adjoint.\\

In \cite{Hoch} was proved that the family
\[\Psi_k(x):=\frac{2^{1/2}}{|J'_\nu(j_{\nu,k})|}x^{1/2}J_\nu(j_{\nu,k}x),\quad k\geq 1,\]
is an orthonormal basis for $L^2(0,1)$.\\ 

Let $\U$ be the unitary operator $\mathcal{U}:L^2(0,1)\rightarrow L^2_\beta(0,1)$ given by
\begin{equation*}
\U u(x):=\kappa_\al^{1/2}x^{-\al/4-\beta/2}u(x^{\kappa_\al}), \quad u\in L^2(0,1).
\end{equation*}
Notice that $\U\Psi_k=\Fi_k$, $k\geq 1$, therefore $\Fi_k$, $k\geq 1$, is an orthonormal basis for $L^2_\beta(0,1)$. We also can see that $\Fi_k\in H^1_{\al,\beta,N}$ by using that $\nu>0,$ (\ref{bessel}) and (\ref{asincero}).\\

Now we set $w(x)=y(z)$ with $z=cx^a$, $a,c>0$. Assume that $y=J_\nu$. Therefore $y$ satisfies the differential equation (\ref{Besselode}), i.e
\[z\frac{d}{dz}\pts{z\frac{dy}{dz}}+(z^2-\nu^2)y=0,\]
which implies that
\[x\frac{d}{dx}\pts{x\frac{dw}{dx}}+a^2(c^2x^{2a}-\nu^2)w=0.\]
Then we set $v(x)=x^bw(x)$, $b\in \R$. Hence
\[x^{2-2a}\frac{d^2v}{dx^2}+(1-2b)x^{1-2a}\frac{dv}{dx}+(b^2-a^2\nu^2 )x^{-2a}v=-a^2c^2v.\]
Finally, we take $a=\kappa_\al, b=(1-\al-\beta)/2$, and $c=j_{\nu,k}$, $k\geq 1$, to get $\Fi_k(1)=0$ and $\A\Phi_k=\lambda_k\Phi_k$ for all $k \geq 1$.
	\end{proof}	

Then $(\A,D(\A))$ is the infinitesimal generator of a diagonalizable analytic semigroup of contractions in $L^2_{\beta}(0,1)$. Thus, we consider interpolation spaces for the initial data. For any $s\geq 0$, we define
\[\mathcal{H}^{s}=\mathcal{H}^{s}(0,1):=D(\mathcal{A}^{s/2})=\left\{u=\sum_{k=1}^\infty a_{k} \Phi_{k}:\|u\|_{\mathcal{H}^{s}}^{2}=\sum_{k=1}^\infty |a_{k}|^{2} \lambda_{k}^{s}<\infty\right\},\]
and we also consider the corresponding dual spaces
\[\mathcal{H}^{-s}:=\left[\mathcal{H}^{s}(0,1)\right]^{\prime}.\]
It is well known that $\mathcal{H}^{-s}$ is the dual space of $\mathcal{H}^{s}$ with respect to the pivot space $L^2_\beta(0,1)$, i.e
\[\mathcal{H}^s\hookrightarrow \mathcal{H}^0=L^2_{\beta}(0,1)=\left(L^2_{\beta}(0,1)\right)'\hookrightarrow \mathcal{H}^{-s},\quad s>0. \]
Equivalently, $\mathcal{H}^{-s}$ is the completion of $L^2_\beta(0,1)$ with respect to the norm
\[\|u\|^2_{-s}:=\sum_{k=1}^{\infty}\lambda_k^{-s}|\langle u,\Phi_k\rangle_\beta|^2.\]
It is well known that the linear mapping given by
\[S(t)u_0=\sum_{k=1}^\infty \textrm{e}^{-\lambda_k t}a_k\Fi_k\quad\text{if}\quad u_0=\sum_{k=1}^\infty a_{k} \Phi_{k}\in \mathcal{H}^s\]
defines a self-adjoint semigroup $S(t)$, $t\geq 0$, in $\mathcal{H}^s $ for all $s\in\R$.\\

For $\delta\in\R$ and a function $z:(0,1)\rightarrow \R$ we introduce the notion of $\delta$-generalized limit of $z$ at $x=0$ as follows
\[\OO_\delta (z):=\lim_{x\rightarrow 0^+} x^\delta z(x).\]

Now we consider a convenient definition of a weak solution for system (\ref{problem}), we multiply the equation in (\ref{problem}) by $x^\beta\varphi(\tau)=x^\beta S(\tau-t)z^{\tau}$, integrate by parts (formally), and take the expression obtained.

\begin{definition}
	Let $T>0$ and $\al,\beta,\mu\in \R$ with $0\leq\al<2$, $\al + \beta >1$, $\mu<\mu(\alpha+\beta)$. Let $f \in L^2(0,T)$ and $u_0\in \HH^{-s}$ for some $s > 0$. A weak solution of (\ref{problem}) is a function $u \in C^0([0,T];\HH^{-s})$ such that for every $\tau \in (0,T]$ and for
	every $z^\tau \in \HH^s$ we have
	\begin{equation}\label{weaksol}
		\left\langle u(\tau), z^{\tau}\right\rangle_{\mathcal{H}^{-s}, \mathcal{H}^{s}}=-\int_{0}^{\tau} f(t) \mathcal{O}_{\al+\beta+\gamma}\left(S(\tau-t) z^{\tau}\right) \mathrm{d} t+\left\langle u_{0}, S(\tau)z^\tau\right\rangle_{\mathcal{H}^{-s}, \mathcal{H}^{s}},
	\end{equation}
	where $\gamma=\gamma(\alpha,\beta,\mu)$ is given in (\ref{gamadef})
\end{definition}
The next result shows the existence of weak solutions for the system (\ref{problem}) under suitable conditions on the parameters $\alpha,\beta,\mu,\gamma$ and $s$. The proof is similar to the proof of Proposition 10 in \cite{GaloLopez}.
\begin{proposition}\label{continuity}
	Let $T>0$ and $\al,\beta\in \R$ with $0\leq\al<2$, $\al + \beta >1$. Let $f \in L^2(0,T)$ and $u_0\in \HH^{-s}$ such that $s>\nu$, where $\nu$ is given in (\ref{Nu}). Then, formula (\ref{weaksol}) defines for each $\tau \in [0, T ]$ a unique element $u(\tau) \in \HH^{-s}$ that can be written as
	\[
	u(\tau)=S(\tau) u_{0}-B(\tau) f, \quad \tau \in(0, T],\]
	where $B(\tau)$ is the strongly continuous family of bounded operators $B(\tau): L^{2}(0,T) \rightarrow \mathcal{H}^{-s}$ given by
	\[\left\langle B(\tau) f, z^{\tau}\right\rangle_{\mathcal{H}^{-s}, \mathcal{H}^{s}}=\int_{0}^{\tau} f(t) \mathcal{O}_{\alpha+\beta+\gamma}\left(S(\tau-t) z^{\tau}\right) \mathrm{d} t, \quad \text{for all  }z^{\tau} \in \mathcal{H}^{s} .\]
	Furthermore, the unique weak solution $u$ on $[0, T]$ to (\ref{problem}) (in the sense of (\ref{weaksol})) belongs to ${C}^{0}\left([0, T] ; \mathcal{H}^{-s}\right)$ and fulfills
	\[
	\|u\|_{L^{\infty}\left([0, T] ; \mathcal{H}^{-s}\right)} \leq C\left(\left\|u_{0}\right\|_{\mathcal{H}^{-s}}+\|f\|_{L^{2}(0, T)}\right).
	\]
\end{proposition}
\begin{proof}
	Fix $\tau>0$. Let $u(\tau)\in \mathcal{H}^{-s}$ be determined by the condition (\ref{weaksol}), hence
	$$
	-u(\tau)+S(\tau) u_{0}=\zeta(\tau)f,
	$$
	where
	$$
	\left\langle\zeta(\tau)f, z^{\tau}\right\rangle_{\mathcal{H}^{-s}, \mathcal{H}^{s}}=\int_{0}^{\tau} f(t) \mathcal{O}_{\alpha+\beta+\gamma}\left(S(\tau-t) z^{\tau}\right) \mathrm{d} t, \quad \text{for all  } z^{\tau} \in \mathcal{H}^{s}.
	$$
	We claim that $\zeta(\tau)$ is a bounded operator from $L^{2}(0, T)$ into $\mathcal{H}^{-s}$: consider $z^{\tau} \in \mathcal{H}^{s}$ given by 
	\begin{equation}\label{finaldata}
		z^{\tau}=\sum_{k=1}^\infty a_{k} \Phi_{k},
	\end{equation}
	therefore
	$$
	S(\tau-t) z^{\tau}=\sum_{k=1}^\infty \mathrm{e}^{\lambda_{k}(t-\tau)} a_{k} \Phi_{k}, \quad \text{for all } t \in[0, \tau].
	$$
	By using Lemma \ref{reduce} and (\ref{lim1}) we obtain that there exists a constant $C=C(\al,\beta,\mu)>0$ such that
	\[|\mathcal{O}_{\alpha+\beta+\gamma}\left(\Phi_k\right)|\leq C |j_{\nu,k}|^{\nu+1/2},\quad k\geq 1,\]
	hence (\ref{below}) implies that there exists a constant $C=C(\al,\beta,\mu)>0$ such that
	\begin{eqnarray*}
		\pts{\int_{0}^{\tau}\left|\mathcal{O}_{\alpha+\beta+\gamma}\left(S(\tau-t) z^{\tau}\right)\right|^{2} \mathrm{~d} t}^{1/2}&\leq & \sum_{k=1}^\infty 
		|a_k| |\OO_{\al+\beta+\gamma}(\Fi_k)| \pts{\int_0^\tau \mathrm{e}^{2\lambda_{k}(t-\tau)} \mathrm{~d} t}^{1/2} \\
		&\leq & C\left\|z^{\tau}\right\|_{\mathcal{H}^{s}}\pts{\sum_{k=1}^\infty |\lambda_k|^{\nu-1/2-s}\pts{1-\mathrm{e}^{-2\lambda_k \tau}}}^{1/2}\\
		&\leq& C\left\|z^{\tau}\right\|_{\mathcal{H}^{s}}\pts{\sum_{k=1}^\infty\frac{1}{k^{2(s-\nu+1/2)}}}^{1/2}=C\left\|z^{\tau}\right\|_{\mathcal{H}^{s}}.
	\end{eqnarray*}
	Therefore $\|\zeta(\tau) f\|_{\mathcal{H}^{-s}}\leq C\|f\|_{L^2(0,T)}$ for all $f\in L^2(0,T)$, $\tau\in (0,T]$.\\
	
	Finally, we fix $f\in L^2(0,T)$ and show that the mapping $\tau\mapsto \zeta(\tau) f$ is right-continuous on $[0,T)$. Let $h>0$ small enough and $z\in \mathcal{H}^s$ given as in (\ref{finaldata}). Thus, proceeding as in the last inequalities, we have
	\begin{eqnarray*}
		|\left\langle\zeta(\tau+h)f-\zeta(\tau)f, z\right\rangle_{\mathcal{H}^{-s}, \mathcal{H}^{s}}|
		&\leq& C\left\|z\right\|_{\mathcal{H}^{s}}\|f\|_{L^2(0,T)}\cts{\pts{\sum_{k=1}^\infty\frac{I(\tau,k,h)}{k^{2(s-\nu+1/2)}}}^{1/2}+\pts{\sum_{k=1}^\infty\frac{1-\mathrm{e}^{-2\lambda_k h}}{k^{2(s-\nu+1/2)}}}^{1/2}},
	\end{eqnarray*}
	where 
	\begin{equation}\label{Icomp}
	I(\tau,k,h)=\lambda_k\int_0^\tau\pts{\mathrm{e}^{\lambda_k(t-\tau-h)}-\mathrm{e}^{\lambda_k(t-\tau)}}^2\mathrm{~d} t
	= \frac{1}{2}(1-\mathrm{e}^{-\lambda_k h})^2(1-\mathrm{e}^{-2\lambda_k \tau})\rightarrow 0\quad\text{as}\quad h\rightarrow 0^+.
	\end{equation}
	Since $0\leq I(\tau,k,h)\leq 1/2$ uniformly for $\tau, h>0$, $k\geq 1$, the result follows by the dominated convergence theorem.
	\end{proof}

\begin{remark}
	In the following section, we will consider initial conditions in $L^2_\beta(0,1)$. Notice that $L^2_\beta(0,1)\subset H^{-\nu-\delta}$ for all $\delta>0$, and we can apply Proposition \ref{continuity} with $s=\nu+\delta$, $\delta>0$, then the corresponding solutions will be in $C^0([0, T ], H^{-\nu-\delta})$.
\end{remark}

\section{Control at the left endpoint}\label{leftend}

\subsection{Upper estimate of the cost of the null controllability}\label{P1Control}
In this section we use the method moment, introduced by Fattorini \& Russell in \cite{Fatorini}, to prove the null controllability of the system (\ref{problem}). In \cite[Section 3]{GaloLopez} the authors construct a biorthogonal family $\d \{\psi_k\}_{k\geq 1}\subset L^2(0,T)$ to the family of exponential functions $\{\mathrm{e}^{-\lambda_{k}(T-t)}\}_{k\geq 1}$ on $[0, T]$, i.e that satisfies
$$
\int_{0}^{T}\psi_k(t) \mathrm{e}^{-\lambda_{l}(T-t)} dt = \delta_{kl},\quad\text{for all}\quad k,l\geq 1.
$$
That construction will help us to get an upper bound for the cost of the null controllability of the system (\ref{problem}). Here, we sketch the process to get the biorthogonal family $\d \lvs{\psi_k}_{k\geq 1}$, see \cite[Section 3]{GaloLopez} for details.\\

Consider the Weierstrass infinite product
\begin{equation}
	\Lambda(z):=\prod_{k=1}^{\infty}\pts{1+\dfrac{iz}{(\kappa_\al j_{\nu, k})^2}}.
\end{equation}
From (\ref{asint}) we have that $j_{\nu, k}=O(k)$ for $k$ large, thus the infinite product is well-defined and converges absolutely in $\C$. Hence $\Lambda(z)$ is an entire function with simple zeros at $i(\kappa_\al j_{\nu, k})^2=i\lambda_k$, $k\geq 1$. It follows that \begin{equation}\label{PsiFunction}
	\Psi_{k}(z):=\dfrac{\Lambda(z)}{\Lambda'(i\lambda_{k})(z-i\lambda_{k})},\quad k\geq 1,
\end{equation}
is an entire function with simple zeros at $i\lambda_k$, $k \geq 1$. Since $\Psi_{k}(x)$ is not in $L^2(\R)$, we need to fix this using a suitable ``complex multiplier", to do this we follow the approach introduced in \cite{Tucsnak}.\\

For $\theta>0$ and $a>0$, we define
$$
\sigma_{\theta}(t):=\exp\pts{-\frac{\theta}{1-t^2}},\quad t\in(-1,1),
$$
and extended by $0$ outside of $(-1, 1)$. Clearly $\sigma_{\theta}$ is analytic on $(-1,1)$. Set $C_{\theta}^{-1}:=\int_{-1}^{1}\sigma_{\theta}(t)\mathrm{d}t$ and define
\begin{equation}\label{Hfunction}
	H_{a,\theta}(z)=C_{\theta}\int_{-1}^{1}\sigma_{\theta}(t)\exp\pts{-iatz}\mathrm{d}t.
\end{equation}
Clearly $H_{a,\theta}(z)$ is an entire function. The following result gives additional information about $H_{a,\theta}(z)$.\\
\begin{lemma}
	The function $H_{a,\theta}$ fulfills the following inequalities
	\begin{eqnarray}
		H_{a,\theta}(ix)&\geq &\frac{\exp\pts{a|x|/\pts{2\sqrt{\theta+1}}}}{11\sqrt{\theta+1}},\quad x\in\R,\label{Hcot1}\\
		|H_{a,\theta}(z)|     &\leq & \exp\pts{a|\Im(z)|},\quad z\in\C,\label{Hcot2}\\
		|H_{a,\theta}(x)|     &\leq & \chi_{|x|\leq 1}(x)+c\sqrt{\theta+1}\sqrt{a\theta\abs{x}}\exp\pts{3\theta/4-\sqrt{a\theta\abs{x}}}\chi_{|x|> 1}(x),\quad x\in\R,\label{Hcot3}
	\end{eqnarray}
	where $c>0$ does not depend on $a$ and $\theta$.
\end{lemma}
We refer to \cite[pp. 85--86]{Tucsnak} for the details.\\

For $k\geq 1$ consider the entire function $F_{k}$ given as
\begin{equation}\label{Ffunction}
	F_{k}(z):=\Psi_{k}(z)\dfrac{H_{a,\theta}(z)}{H_{a,\theta}(i\lambda_{k})},\quad z\in\C.
\end{equation}

For $\delta\in(0,1)$ we set
\begin{equation}\label{aConst}
	a:=\frac{T(1-\delta)}{2}>0,\quad \text{and}\quad \theta:=\dfrac{(1+\delta)^2}{\kappa_\al^2 T\pts{1-\delta}}>0.
\end{equation}

\begin{lemma}
	For each $k\geq 1$ the function $F_{k}(z)$ satisfies the following properties:\\
	i) $F_{k}$ is of exponential type $T/2$.\\
	ii) $F_{k}\in L^1(\R)\cap L^2(\R)$.\\
	iii) $F_k$ satisfies $F_{k}(i\lambda_{l})=\delta_{kl}$ for all $k,l\geq 1$.\\
	iv) Furthermore, there exists a constant $c>0$, independent of $T,\alpha$ and $\delta$, such that
	\begin{equation}\label{Fbound}
		\left\|F_{k}\right\|_{L^{1}(\R)} \leq \frac{C(T, \alpha,\delta)}{\lambda_{k}\left|\Lambda^{\prime}\left(i \lambda_{k}\right)\right|}  \exp\pts{-\frac{a\lambda_k}{2\sqrt{\theta+1}}},		
	\end{equation} 
	where
	\begin{equation}\label{upper}
		C(T, \alpha,\delta)=c\sqrt{\theta+1}\cts{\exp\pts{{\frac{1}{\sqrt{2}\kappa_\alpha}}}+\sqrt{\theta+1}\frac{\kappa_\alpha^2}{\delta^3}\exp\pts{\frac{3 \theta}{4}}}. 
	\end{equation}
\end{lemma}

The $L^2$-version of the Paley-Wiener theorem implies that there exists $\eta_k\in L^2(\R)$ with support in $[-T/2,T/2]$ such that $F_k(z)$ is the analytic extension of the Fourier transform of $\eta_k$. Hence 
\begin{equation}\label{psieta}
	\psi_k(t):=\mathrm{e}^{\lambda_k T/2}\eta_k(t-T/2),\quad t\in[0,T],\,k\geq1,
\end{equation}
is the family we are looking for.\\

Since $\eta_k, F_k\in L^1(\R)$, the inverse Fourier theorem yields 
\[\eta_k(t)=\frac{1}{2\pi}\int_{\R}\mathrm{e}^{it\tau}F_k(\tau)\mathrm{d}\tau,\quad t\in\R, k\geq 1,\]
hence (\ref{psieta}) implies that $\psi_k\in C([0,T])$, and by using (\ref{Fbound}) we have
\begin{equation}\label{psiL1}
	\|\psi_k\|_{\infty}\leq \frac{C(T, \alpha,\delta)}{\lambda_{k}\left|\Lambda^{\prime}\left(i \lambda_{k}\right)\right|}  \exp\pts{\frac{T\lambda_k}{2}-\frac{a\lambda_k}{2\sqrt{\theta+1}}},\quad k\geq 1.	
\end{equation}

Now, we are ready to prove the null controllability of the system (\ref{problem}). Let $u_{0}\in L^{2}_\beta(0,1)$. Then consider its Fourier series with respect to the orthonormal basis $\{\Phi_{k}\}_{k\geq 1}$,
\begin{equation}\label{uoSerie}
	u_{0}(x)=\sum_{k=1}^{\infty} a_{k} \Phi_{k}(x).
\end{equation}
We set
\begin{equation}\label{fserie}
	f(t):=\sum_{k=1}^{\infty}\frac{a_{k} \mathrm{e}^{-\lambda_{k} T}}{\mathcal{O}_{\alpha+\beta+\gamma}\left(\Phi_{k}\right)} \psi_{k}(t).
\end{equation}
Since $\{\psi_k\}$ is biorthogonal to $\{\mathrm{e}^{-\lambda_k(T-t)}\}$ we have
$$
\int_{0}^{T} f(t) \mathcal{O}_{\alpha+\beta+\gamma}\left(\Phi_{k}\right) \mathrm{e}^{-\lambda_{k}(T-t)} \mathrm{d} t = a_{k} \mathrm{e}^{-\lambda_{k} T} = \left\langle u_{0}, \mathrm{e}^{-\lambda_{k} T}\Phi_{k}\right\rangle_\beta = \left\langle u_{0}, \mathrm{e}^{-\lambda_{k} T}\Phi_{k}\right\rangle_{\mathcal{H}^{-s}, \mathcal{H}^{s}}.
$$
Let $u\in C([0,T];H^{-s})$ that satisfies (\ref{weaksol}) for all $\tau\in (0,T]$, $z^\tau\in H^s$. In particular, for $\tau=T$ we take $z^T=\Phi_k$, $k\geq 1$, then the last equality implies that
$$
\left\langle u(\cdot, T), \Phi_{k}\right\rangle_{\mathcal{H}^{-s}, \mathcal{H}^{s}}=0\quad \text{for all}\quad k \geq 1,
$$
hence $u(\cdot, T)\equiv 0$.\\

It just remains to estimate the norm of the control $f$. From (\ref{psiL1}) and (\ref{fserie})  we get
\begin{equation}\label{finty}
	\|f\|_{\infty} \leq  C(T, \alpha,\delta)\sum_{k=1}^{\infty} \frac{\left|a_{k}\right|}{\left|\mathcal{O}_{\alpha+\beta+\gamma}\left(\Phi_{k}\right)\right|} \frac{1}{\lambda_{k}\left|\Lambda^{\prime}\left(i \lambda_{k}\right)\right|} \exp\pts{-\frac{T\lambda_k}{2}-\frac{a\lambda_k}{2\sqrt{\theta+1}}}.
\end{equation}

Using \cite[Chap. XV, p. 438, eq. (3)]{Watson}, we can write
\begin{equation*}
	\Lambda(z)=\Gamma(\nu+1)\pts{\dfrac{2\kappa_\al}{\sqrt{-iz}}}^\nu J_{\nu}\pts{\dfrac{\sqrt{-iz}}{\kappa_\al}},
\end{equation*}
therefore
\begin{equation}\label{deriLambda}
	\left|\Lambda^{\prime}\left(i \lambda_{k}\right)\right|=\Gamma(\nu+1)\frac{2^{\nu}}{|j_{\nu, k}|^{\nu}} \frac{1}{2\kappa_\al^{2} j_{\nu, k}} |J_{\nu}^{\prime}\left(j_{\nu, k}\right)|, \quad k\geq 1,
\end{equation}
and by using (\ref{lambdak}) and (\ref{lim1}) we get
\begin{equation*}\label{estimate2}
	\d\left|\mathcal{O}_{\alpha+\beta+\gamma}\left(\Phi_{k}\right)\lambda_{k}\Lambda^{\prime}\left(i \lambda_{k}\right)\right|= 2^{-1/2}\sqrt{\kappa_\al} j_{\nu, k}. 
\end{equation*}
From (\ref{finty}), the last two equalities and using that $\lambda_k\geq \lambda_1$, it follows that 
\begin{equation*}
	\|f\|_{\infty} \leq \frac{C(T, \alpha,\delta)}{\sqrt{\kappa_\al}} \exp\pts{-\frac{T\lambda_1}{2}-\frac{a\lambda_1}{2\sqrt{\theta+1}}}\sum_{k=1}^{\infty} \frac{|a_{k}| }{j_{\nu, k}}.
\end{equation*}
By using the Cauchy-Schwarz inequality, the fact that $j_{\nu,k}\geq (k-1/4)\pi$ (by (\ref{below}))  and (\ref{uoSerie}), we obtain that
\begin{eqnarray*}
	\|f\|_{\infty} &\leq& \frac{C(T, \alpha,\delta)}{\sqrt{\kappa_\al}} \exp\pts{-\frac{T\lambda_1}{2}-\frac{a\lambda_1}{2\sqrt{\theta+1}}}\left\|u_{0}\right\|_{\beta}.
\end{eqnarray*}
Notice that $0<\kappa_\alpha\leq 1$, and $\theta>0$. Thus, by using (\ref{aConst}) with $\delta\in(0,1)$, we have that
\[\theta\leq \frac{4}{(1-\delta)\kappa_\alpha^2 T },\quad \sqrt{\theta+1}\leq \frac{2(1+T)^{1/2}}{(1-\delta)^{1/2}\kappa_\alpha T^{1/2}},\quad \sqrt{\theta+1}\leq \theta+1,\]
therefore
\begin{equation}\label{reduc}
\frac{a}{\sqrt{\theta+1}}\geq \frac{\kappa_\alpha(1-\delta)^{3/2}T^{3/2}}{4(1+T)^{1/2}}, \quad C(T, \alpha,\delta)\leq c\pts{1+\frac{1}{(1-\delta)\kappa_\alpha^2 T}}\cts{\exp\pts{\frac{1}{\sqrt{2}\kappa_\alpha}}+\frac{1}{\delta^3}\exp\pts{\frac{3}{(1-\delta)\kappa_\alpha^2 T}}},
\end{equation}
and by using the definition of $\lambda_{1}$ the result follows.

\subsection{Lower estimate of the cost of the null controllability}
In this section, we get a lower estimate of the cost $\mathcal{K}=\mathcal{K}(T,\alpha,\beta,\mu)$.
We set
\begin{equation}\label{first} 
	u_0(x):= \frac{\abs{J'_\nu\pts{j_{\nu,1}}}}{(2\kappa_\al)^{1/2}}\Fi_1(x),\,x\in(0,1),
	\quad \text{hence}\quad\|u_0\|^2_\beta=\frac{\abs{J'_\nu\pts{j_{\nu,1}}}^2}{2\kappa_\al}.
\end{equation}
For $\varepsilon>0$ small enough, there exists $f\in U(\al,\beta,\mu,T,u_0)$ such that
\begin{equation}\label{inicost}
	u(\cdot,T)\equiv 0,\quad \text{and}\quad \|f\|_{L^2(0,T)}\leq (\mathcal{K}+\varepsilon)\|u_0\|_\beta.
\end{equation}
Then, in (\ref{weaksol}) we set $\tau=T$ and  take $z^\tau=\Fi_k$, $k\geq 1$,  to obtain
\begin{eqnarray*}
	\mathrm{e}^{-\lambda_k T}\left\langle u_{0},\Fi_k\right\rangle_\beta=\left\langle u_{0}, S(T)\Fi_k\right\rangle_{\mathcal{H}^{-s}, \mathcal{H}^{s}}&=&
	\int_{0}^{T} f(t) \mathcal{O}_{\al+\beta+\gamma}\left(S(T-t) \Fi_k\right) \mathrm{d} t\\
	&=&
	\mathrm{e}^{-\lambda_k T}\mathcal{O}_{\al+\beta+\gamma}\left(\Fi_k\right)\int_{0}^{T} f(t) \mathrm{e}^{\lambda_k t} \mathrm{d} t,
\end{eqnarray*}
from (\ref{first}) and (\ref{lim1}) it follows that
\begin{equation}\label{orto}
	\int_{0}^{T} f(t) \mathrm{e}^{\lambda_k t} \mathrm{d} t= \frac{2^{\nu}\Gamma(\nu+1)\abs{J'_{\nu}\pts{j_{\nu, 1}}}^2}{2\kappa_\al\pts{j_{\nu,1}}^{\nu}}\delta_{1,k},\quad k\geq 1.
\end{equation}

Now consider the function $v: \mathbb{C} \rightarrow \mathbb{C}$ given by
\begin{equation}\label{vanalitica}
	v(s):=\int_{-T / 2}^{T / 2} f\left(t+\frac{T}{2}\right) \mathrm{e}^{-i s t} \mathrm{~d} t, \quad s \in \mathbb{C} .
\end{equation}
Fubini and Morera's theorems imply that $v(s)$ is an entire function. Moreover, (\ref{orto}) implies that
\[v(i\lambda_k)=0\quad\text{for all }k\geq 2,\quad \text{and}\quad v(i\lambda_1)=\frac{2^{\nu}\Gamma(\nu+1)\abs{J'_{\nu}\pts{j_{\nu, 1}}}^2}{2\kappa_\al\pts{j_{\nu,1}}^{\nu}}\mathrm{e}^{-\lambda_1 T/2}.\]
We also have that
\begin{equation}\label{uve}
	|v(s)| \leq \mathrm{e}^{T|\Im(s)|/2}\int_{0}^{T}|f(t)| \mathrm{d} t \leq (\mathcal{K}+\varepsilon) T^{1/2}\mathrm{e}^{T|\Im(s)|/2} \left\|u_{0}\right\|_{\beta}.
\end{equation}
Consider the entire function $F(z)$ given by
\begin{equation}\label{entire}
	F(s):=v\left(s-i \delta\right), \quad s \in \mathbb{C},
\end{equation}
for some $\delta>0$ that will be chosen later on. Clearly, 
\[ F\left(a_{k}\right)=0, \quad k\geq 2, \quad \text {where} \quad a_{k}:=i\left(\lambda_k+\delta\right),\quad k\geq 1,\quad\text{and}\]
\begin{equation}\label{ena1}
	F\left(a_{1}\right)= \frac{2^{\nu}\Gamma(\nu+1)\abs{J'_{\nu}\pts{j_{\nu, 1}}}^2}{2\kappa_\al\pts{j_{\nu,1}}^{\nu}}\mathrm{e}^{-\lambda_1 T/2}.
\end{equation}
From (\ref{first}), (\ref{uve}) and (\ref{entire}) we obtain
\begin{equation}\label{logF}
	\log |F(s)|\leq \frac{T}{2}|\Im(s)-\delta|+\log\pts{(\mathcal{K}+\varepsilon) T^{1 / 2}\frac{\abs{J'_\nu\pts{j_{\nu,1}}}}{\pts{2\kappa_\al}^{1/2}}},\quad s\in\mathbb{C}.
\end{equation}

We apply Theorem \ref{repreTheo} to the function $F(z)$ given in (\ref{entire}). In this case, (\ref{uve}) implies that $A\leq T/2$. Also notice that $\Im\left(a_{k}\right)>0$, $k\geq 1$, to get
\begin{equation}\label{aprep}
	\log \left|F\left(a_{1}\right)\right|\leq\left(\lambda_1+\delta\right)\frac{T}{2}+\sum_{k=2}^{\infty} \log \left|\frac{a_{1}-a_{k}}{a_{1}-\bar{a}_{k}}\right|+\frac{\Im\left(a_{1}\right)}{\pi} \int_{-\infty}^{\infty} \frac{\log |F(s)|}{\left|s-a_{1}\right|^{2}} \mathrm{~d}s.
\end{equation}
By using the definition of the constants $a_k$'s we have
\begin{eqnarray}
	\sum_{k=2}^{\infty} \log \left|\frac{a_{1}-a_{k}}{a_{1}-\bar{a}_{k}}\right|&=&\sum_{k=2}^{\infty} \log \left(\frac{\left( j_{\nu, k}\right)^{2}-\left( j_{\nu, 1}\right)^{2}}{2 \delta / \kappa_\alpha^{2}+\left( j_{\nu, 1}\right)^{2}+\left(j_{\nu, k}\right)^{2}}\right)\notag\\
	&\leq& \sum_{k=2}^{\infty} \frac{1}{j_{\nu, k+1}-j_{\nu, k}} \int_{j_{\nu, k}}^{j_{\nu, k+1}} \log \left(\frac{ x^{2}}{2 \delta / \kappa_\alpha^{2}+ x^{2}}\right) \mathrm{d} x \label{apoyo1} \\ 
	&\leq& \frac{1}{\pi} \int_{j_{\nu, 2}}^{\infty} \log \left(\frac{ x^{2}}{2 \delta / \kappa_\alpha^{2}+x^{2}}\right) \mathrm{d} x,\notag\\
	&=& -\frac{j_{\nu, 2}}{\pi} \log \left(\frac{1}{1+2 \delta /\left(\kappa_\alpha j_{\nu, 2}\right)^{2}}\right)- \frac{2\sqrt{2 \delta}}{\pi\kappa_\alpha}\left(\frac{\pi}{2}-\tan ^{-1}\left(\kappa_\alpha j_{\nu, 2} / \sqrt{2 \delta}\right)\right) ,\notag
\end{eqnarray}
where we have used Lemma \ref{consec} and made the change of variables
$$
\tau=\frac{ \kappa_\alpha}{\sqrt{2 \delta}} x.
$$
From (\ref{logF}) we get the estimate
\begin{equation}\label{apoyo2}
	\frac{\Im\left(a_{1}\right)}{\pi} \int_{-\infty}^{\infty} \frac{\log |F(s)|}{\left|s-a_{1}\right|^{2}} \mathrm{~d} s \leq \frac{T \delta}{2}+\log \left((\mathcal{K}+\varepsilon) T^{1 / 2} \frac{\left|J_{\nu}^{\prime}\left(j_{\nu, 1}\right)\right|}{\pts{2 \kappa_\alpha}^{1/2}}\right).
\end{equation}
From (\ref{ena1}), (\ref{aprep}), (\ref{apoyo1}) and (\ref{apoyo2}) we have
\begin{equation}\label{abajo}
\frac{2\sqrt{2 \delta}}{\pi\kappa_\alpha}\tan ^{-1}\left(\frac{\sqrt{2 \delta}}{\kappa_\alpha j_{\nu, 2}}\right) -\frac{j_{\nu, 2}}{\pi} \log \left(1+\frac{2 \delta}{ \left(\kappa_\alpha j_{\nu, 2}\right)^{2}}\right) 
-\pts{\lambda_1+\delta}T \leq \log(\mathcal{K}+\varepsilon)+\log h(\alpha,\beta,\mu, T),
\end{equation}
where
\[h(\alpha,\beta,\mu, T)=\frac{\pts{{2T \kappa_\alpha}}^{1/2}\left(j_{\nu, 1}\right)^{\nu}}{2^{\nu} \Gamma(\nu+1) \left|J_{\nu}^{\prime}\left(j_{\nu, 1}\right)\right|}.\]
The result follows by taking 
\[\delta=\frac{\kappa_\alpha^2 \pts{j_{\nu,2}}^2}{2}, \quad \text{and then letting}\quad \varepsilon\rightarrow 0^+.\]

\section{Control at the right endpoint}\label{rightend}
Here, we analyze the null controllability of the system (\ref{problem2}) where $\alpha +\beta >1$, $0\leq \al < 2$,  $\mu$ and $\gamma$ satisfy (\ref{mucon}) and (\ref{gamadef}) respectively. As in Section \ref{leftend} we give a suitable definition of a weak solution for the system (\ref{problem2}).\\

\begin{definition}
	Let $T>0$ and $\al,\beta,\mu\in \R$ with $0\leq\al<2$, $\al + \beta >1$, $\mu<\mu(\alpha+\beta)$. Let $f \in L^2(0,T)$ and $u_0\in \HH^{-s}$ for some $s > 0$. A weak solution of (\ref{problem2}) is a function $u \in C^0([0,T];\HH^{-s})$ such that for every $\tau \in (0,T]$ and for
	every $z^\tau \in \HH^s$ we have
	\begin{equation}\label{weaksol2}
		\left\langle u(\tau), z^{\tau}\right\rangle_{\mathcal{H}^{-s}, \mathcal{H}^{s}}=\left\langle u_{0}, S(\tau)z^\tau\right\rangle_{\mathcal{H}^{-s}, \mathcal{H}^{s}} - \int_{0}^{\tau} f(t)\lim_{x\rightarrow 1^- }S(\tau-t) z^{\tau}_{x}(x) \mathrm{d} t,
	\end{equation}
	where $\gamma=\gamma(\alpha,\beta,\mu)$ is given by (\ref{gamadef}).
\end{definition}
The next result shows the existence of weak solutions for the system (\ref{problem2}) under certain conditions on the parameters $\alpha,\beta,\mu,\gamma$ and $s$.
\begin{proposition}\label{continuity2}
	Let $T>0$ and $\al,\beta\in \R$ with $0\leq\al<2$, $\al + \beta >1$. Let $f \in L^2(0,T)$ and $u_0\in \HH^{-s}$ such that $s>1/2$. Then, formula (\ref{weaksol2}) defines for each $\tau \in [0, T ]$ a unique element $u(\tau) \in \HH^{-s}$ that can be written as
	\[
	u(\tau)=S(\tau) u_{0}-B(\tau) f, \quad \tau \in(0, T],\]
	where $B(\tau)$ is the strongly continuous family of bounded operators $B(\tau): L^{2}(0,T) \rightarrow \mathcal{H}^{-s}$ given by
	\[\left\langle B(\tau) f, z^{\tau}\right\rangle_{\mathcal{H}^{-s}, \mathcal{H}^{s}}=\int_{0}^{\tau} f(t)\lim_{x\rightarrow 1^- }S(\tau-t) z^{\tau}_{x}(x) \mathrm{d} t, \quad \text{for all  }z^{\tau} \in \mathcal{H}^{s} .\]
	Furthermore, the unique weak solution $u$ on $[0, T]$ to (\ref{problem2}) (in the sense of (\ref{weaksol2})) belongs to ${C}^{0}\left([0, T] ; \mathcal{H}^{-s}\right)$ and fulfills
	\[
	\|u\|_{L^{\infty}\left([0, T] ; \mathcal{H}^{-s}\right)} \leq C\left(\left\|u_{0}\right\|_{\mathcal{H}^{-s}}+\|f\|_{L^{2}(0, T)}\right).
	\]
\end{proposition}
\begin{proof}
	Fix $\tau>0$. Let $u(\tau)\in H^{-s}$ be determined by the condition (\ref{weaksol2}), hence
	$$
	-u(\tau)+S(\tau) u_{0}=\zeta(\tau)f,
	$$
	where
	$$
	\left\langle\zeta(\tau)f, z^{\tau}\right\rangle_{\mathcal{H}^{-s}, \mathcal{H}^{s}}=\int_{0}^{\tau} f(t) \lim_{x\rightarrow 1^- }S(\tau-t) z^{\tau}_{x}(x)\mathrm{d} t \quad \text{for all  } z^{\tau} \in \mathcal{H}^{s}.
	$$
	Let $z^{\tau} \in \mathcal{H}^{s}$ given by 
	\begin{equation}\label{finaldata2}
		z^{\tau}=\sum_{k=1}^\infty a_{k} \Phi_{k},
	\end{equation}
	therefore
	$$
	\lim_{x\rightarrow 1^- }S(\tau-t) z^{\tau}_{x}(x)=\sum_{k=1}^\infty \mathrm{e}^{\lambda_{k}(t-\tau)} a_{k} \Phi'_{k}(1) \quad \text{for all } t \in[0, \tau].
	$$
	By (\ref{Phik}) we get
	\begin{equation}\label{Dphik}
		\left|\Phi'_{k}(1)\right| = 2^{1/2}\kappa_{\al}^{3/2}j_{\nu, k}, \quad k\geq 1,
	\end{equation}
	hence (\ref{below}) implies that there exists $C=C(\al,\beta,\mu)>0$ such that
	\begin{eqnarray*}
		\pts{\int_{0}^{\tau}\left|\lim_{x\rightarrow 1^- }S(\tau-t) z^{\tau}_{x}(x)\right|^{2} \mathrm{~d} t}^{1/2}&\leq & \sum_{k=1}^\infty 
		|a_k| |\Phi'_{k}(1)| \pts{\int_0^\tau \mathrm{e}^{2\lambda_{k}(t-\tau)} \mathrm{~d} t}^{1/2} \\
		&\leq& C\left\|z^{\tau}\right\|_{\mathcal{H}^{s}}\pts{\sum_{k=1}^\infty |\lambda_k|^{1-s}\int_0^\tau \mathrm{e}^{2\lambda_{k}(t-\tau)} \mathrm{~d} t}^{1/2}\\
		&\leq& C\left\|z^{\tau}\right\|_{\mathcal{H}^{s}}\pts{\sum_{k=1}^\infty\frac{1}{k^{2s}}}^{1/2}=C\left\|z^{\tau}\right\|_{\mathcal{H}^{s}}.
	\end{eqnarray*}
	Therefore $\|\zeta(\tau) f\|_{\mathcal{H}^{-s}}\leq C\|f\|_{L^2(0,T)}$ for all $f\in L^2(0,T)$, $\tau\in (0,T]$.\\
	
	Finally, we fix $f\in L^2(0,T)$ and show that the mapping $\tau\mapsto \zeta(\tau) f$ is right-continuous on $[0,T)$. Let $h>0$ small enough and $z\in \mathcal{H}^s$ given as in (\ref{finaldata2}). Thus, proceeding as in the last inequalities, we have
	\begin{eqnarray*}
		|\left\langle\zeta(\tau+h)f-\zeta(\tau)f, z\right\rangle_{\mathcal{H}^{-s}, \mathcal{H}^{s}}|&\leq&\int_{0}^{\tau} |f(t)|\left|\lim_{x\rightarrow 1^- }(S(\tau+h-t)-S(\tau-t)) z^{\tau}_{x}(x)\right| \mathrm{d} t \\
		&&+ \int_{\tau}^{\tau+h} |f(t)|\left|\lim_{x\rightarrow 1^- }S(\tau+h-t) z^{\tau}_{x}(x)\right| \mathrm{d} t\\
		&\leq& C\left\|z\right\|_{\mathcal{H}^{s}}\|f\|_{L^2(0,T)}\cts{\pts{\sum_{k=1}^\infty\frac{I(\tau,k,h)}{k^{2s}}}^{1/2}+\pts{\sum_{k=1}^\infty\frac{1-\mathrm{e}^{-2\lambda_k h}}{k^{2s}}}^{1/2}},
	\end{eqnarray*}
	where $I(\tau,k,h)$ satisfies (\ref{Icomp}). 
\end{proof}

\begin{remark}
	In the following subsections, we will consider initial conditions in $L^2_\beta(0,1)$. We can apply Proposition \ref{continuity2} with $s=1/2+\delta$, $\delta>0$, then the corresponding solutions will be in $C^0([0, T ], H^{-1/2-\delta})$.
\end{remark}

\subsection{Upper estimate of the cost of the null controllability}
We are ready to prove the null controllability of the system (\ref{problem2}). Let $u_{0}\in L^{2}_\beta(0,1)$ given as follows
\begin{equation}\label{uoSerie2}
	u_{0}(x)=\sum_{k=1}^{\infty} a_{k} \Phi_{k}(x).
\end{equation}

We set
\begin{equation}\label{fserie2}
	f(t):=\sum_{k=1}^{\infty}\frac{a_{k} \mathrm{e}^{-\lambda_{k} T}}{\Phi'_{k}(1)} \psi_{k}(t).
\end{equation}

 Since the sequence $\{\psi_k\}$ is biorthogonal to $\{\mathrm{e}^{-\lambda_k(T-t)}\}$ we have
\begin{equation}\label{control2}
\Phi'_{k}(1)\int_{0}^{T} f(t)\mathrm{e}^{-\lambda_{k}(T-t)} \mathrm{d} t = a_{k} \mathrm{e}^{-\lambda_{k} T} = \left\langle u_{0}, \mathrm{e}^{-\lambda_{k} T}\Phi_{k}\right\rangle_\beta = \left\langle u_{0}, \mathrm{e}^{-\lambda_{k} T}\Phi_{k}\right\rangle_{\mathcal{H}^{-s}, \mathcal{H}^{s}}.
\end{equation}

Let $u\in C([0,T];H^{-s})$ be the weak solution of system (\ref{problem2}).  In particular, for $\tau=T$ we take $z^T=\Phi_k$, $k\geq 1$, then (\ref{weaksol2}) and (\ref{control2}) imply that 
$$
\left\langle u(\cdot, T), \Phi^{k}\right\rangle_{\mathcal{H}^{-s}, \mathcal{H}^{s}}=0\quad \text{for all}\quad k \geq 1,
$$
therefore $u(\cdot, T) = 0$.\\

It just remains to estimate the norm of the control $f$. From (\ref{psiL1}), (\ref{deriLambda}), (\ref{Dphik}) and (\ref{fserie2})  we get
\begin{equation*}
	\|f\|_{\infty} \leq \frac{C(T, \alpha,\delta)\kappa_\al^{1/2}}{2^{\nu}\Gamma(\nu+1)}\sum_{k=1}^{\infty} \frac{|j_{\nu,k}|^{\nu}}{|J'_{\nu}(j_{\nu,k})|}\dfrac{|a_{k}|}{\lambda_{k}}\exp\pts{-\frac{T\lambda_k}{2}-\frac{a\lambda_k}{2\sqrt{\theta+1}}}.
\end{equation*}
By using that $\mathrm{e}^{-x}\leq \mathrm{e}^{-r}r^rx^{-r}$ for all $x,r>0,$ the Cauchy-Schwarz inequality, Lemma \ref{reduce} and the fact that $j_{\nu,k}\geq (k-1/4)\pi$ (by (\ref{below}))  and (\ref{uoSerie}), we obtain that
\begin{eqnarray*}
	\|f\|_{\infty} &\leq& \frac{C(T, \alpha,\delta)}{\pts{2\kappa_\al}^{\nu}\Gamma(\nu+1)}\pts{\dfrac{2\nu+1}{T}}^{(2\nu+1)/4} \exp\pts{-\frac{2\nu+1}{4}}\exp\pts{-\frac{a\lambda_1}{2\sqrt{\theta+1}}-\frac{T\lambda_1}{4}}\sum_{k=0}^{\infty}\frac{|a_k|}{\lambda_k}\\
	&\leq& \frac{C(T, \alpha,\delta)}{\pts{2\kappa_\al}^{\nu}\Gamma(\nu+1)}\pts{\dfrac{2\nu+1}{T}}^{(2\nu+1)/4} \exp\pts{-\frac{2\nu+1}{4}}\exp\pts{-\frac{a\lambda_1}{2\sqrt{\theta+1}}-\frac{T\lambda_1}{4}}\left\|u_{0}\right\|_{\beta},
\end{eqnarray*}
and the result follows by (\ref{reduc}).
\subsection{Lower estimate of the cost of the null controllability at $x=1$}
Here, we just give a sketch of the proof of a lower estimate for the cost  $\widetilde{\mathcal{K}}=\widetilde{\mathcal{K}}(T,\alpha,\beta,\mu)$. Consider $u_0\in L^2_\beta(0,1)$ given in (\ref{first}).\\

For $\varepsilon>0$ small enough, there exists $f\in \widetilde{U}(\al,\beta,\mu,T,u_0)$ such that
\begin{equation}\label{inicost2}
	u(\cdot,T)\equiv 0,\quad \text{and}\quad \|f\|_{L^2(0,T)}\leq (\widetilde{\mathcal{K}}+\varepsilon)\|u_0\|_\beta.
\end{equation}
Then, in (\ref{weaksol2}) we set $\tau=T$ and  take $z^T=\Fi_k$, $k\geq 1$,  to obtain
\[ \mathrm{e}^{-\lambda_k T}\left\langle u_{0},\Fi_k\right\rangle_\beta=\left\langle u_{0}, S(T)\Fi_k\right\rangle_{\mathcal{H}^{-s}, \mathcal{H}^{s}} = \mathrm{e}^{-\lambda_k T}\Fi'_k(1)\int_{0}^{T} f(t) \mathrm{e}^{\lambda_k t} \mathrm{d} t,\]
from (\ref{first}) and (\ref{Dphik}) it follows that
\begin{equation}\label{orto2}
	\int_{0}^{T} f(t) \mathrm{e}^{\lambda_k t} \mathrm{d} t = \frac{\abs{J'_{\nu}\pts{j_{\nu, 1}}}}{2\kappa^{2}_{\al}j_{\nu,1}}\delta_{1,k},\quad k\geq 1.
\end{equation}

Next, we proceed as in (\ref{vanalitica})--(\ref{abajo}). But in this case, the corresponding functions $v$ and $F$ satisfy
\[v(i\lambda_k)=0\quad\text{for all }k\geq 2,\quad v(i\lambda_1)=\frac{\abs{J'_{\nu}\pts{j_{\nu, 1}}}}{2\kappa^{2}_{\al}j_{\nu,1}}\mathrm{e}^{-\lambda_1 T/2},\quad \text{and}\]
\[F(a_k)=0\quad\text{for all }k\geq 2,\quad F(a_1)=\frac{\abs{J'_{\nu}\pts{j_{\nu, 1}}}}{2\kappa^{2}_{\al}j_{\nu,1}}\mathrm{e}^{-\lambda_1 T/2}.\]

Hence we can see that
\[ \frac{2\sqrt{2 \delta}}{\pi\kappa_\alpha}\tan ^{-1}\left(\frac{\sqrt{2 \delta}}{\kappa_\alpha j_{\nu, 2}}\right) -\frac{j_{\nu, 2}}{\pi} \log \left(1+\frac{2 \delta}{ \left(\kappa_\alpha j_{\nu, 2}\right)^{2}}\right) 
-\pts{\lambda_1+\delta}T \leq \log(\widetilde{\mathcal{K}}+\varepsilon)+\log \widetilde{h}(\alpha,\beta,\mu, T),\]
where $\widetilde{h}(\alpha,\beta,\mu, T)=T^{1/2}\kappa^{3/2}_{\al}j_{\nu,1}/\sqrt{2}$.
The result follows by taking $\delta=\kappa_\alpha^2 \pts{j_{\nu,2}}^2/2$ and then letting $\varepsilon\rightarrow 0^+$.

\section{The case $\al+\beta=1$}\label{sumauno}
Concerning the case $\alpha+\beta<1$, in \cite{GaloLopez} we showed the system (\ref{problem}) is well-posed when considering suitable weighted Dirichlet condition at the left endpoint and proved the null-controllability of the corresponding system. In both cases ($\alpha+\beta <1$ and $\alpha+\beta >1$) our approach is based on the validity of the Hardy inequality, see Proposition \ref{Hardy} and \cite[Proposition 4]{GaloLopez}. If $\alpha+\beta =1$ then $\mu(\alpha+\beta)=0$, and the corresponding Hardy inequality does not provide any information. Thus, to solve the case $\al+\beta=1$ we use the singular Sturm-Liouville theory, see \cite{Zettl} for the definitions used here.
\subsection{Singular Sturm-Liouville theory}
Assume that $0\leq \al<2$, and $\mu < 0$. Consider the differential expression $M$ defined by
\[Mu=-(pu')'+qu\]
where $\d p(x) = x, q(x) = -\mu x^{-1}$, and $w(x) = x^{1-\al}$.\\

Clearly,
\begin{equation*}\label{Acond1}
1/p, q, w \in L_{\text{loc}}(0,1),\quad p,w >0\text{ on } (0,1),
\end{equation*} 
thus $Mu$ is defined a.e. for functions $u$ such that $u, pu'\in AC_{\text{loc}}(0,1)$, where $AC_{\text{loc}}(0,1)$ is the space of all locally absolutely continuous functions in $(0,1)$.\\

When $\beta=1-\al$ the operator $\mathcal{A}$ given in (\ref{Aoperator}) can be written as $\A =w^{-1}M$. Now, consider
\begin{equation*}\label{Dmax}
	D_{\max}:=\left\{u\in AC_{\text{loc}}(0,1)\, |\,pu'\in AC_{\text{loc}}(0,1),\, u, \A u\in L^2_{1-\al}(0,1)\right\},\quad\text{and}
\end{equation*}
\[D(\A):=\left\{\begin{aligned}
\{u\in D_{\max}\,| \lim_{x\rightarrow 0^+}x^{\sqrt{-\mu}}u(x)=0,u(1)=0\} & &\text{if } \sqrt{-\mu}<\kappa_\al, \\
 \{u\in D_{\max}\,| u(1)=0\}& & \text{if } \sqrt{-\mu}\geq\kappa_\al.
\end{aligned}\right.
\]
Recall that the Lagrange form is given as follows
\[[u,v]:=upv'-u'pv,\quad u,v\in D_{\max}.\]
\begin{proposition}
Let $0\leq \al<2$, $\mu < 0$, and $\nu=\sqrt{-\mu}/\kappa_\al$. Then $\A:D(\A)\subset L^2_{1-\al}(0,1)\rightarrow L^2_{1-\al}(0,1)$ is a self-adjoint operator. Furthermore, the family given in (\ref{Phik}) is an orthonormal basis for $L^2_{1-\al}(0,1)$ such that
	\begin{equation*}
		\mathcal{A} \Fi_k=\lambda_k \Fi_k, \quad k\geq 1.
	\end{equation*}
\end{proposition}
\begin{proof}
First, we refer to \cite[Definition 7.3.1]{Zettl}.\\ 
Since $1/p,q,w\in L^1(1/2,1)$ we have that $x=1$ is a regular point.
Consider the following functions
\begin{equation*}\label{PrincipalSol0}
y_{+}(x) = x^{\sqrt{-\mu}},\quad y_{-}(x) = x^{-\sqrt{-\mu}}.
\end{equation*} 
Notice that $My_{\pm}=0y_{\pm}$. Since $y_{\pm}>0$ on $(0,1)$ we have that $x=0$ is non-oscillatory (NO) for $\lambda=0$. Theorem 2.2 in \cite{Niessen} implies that $y_+$ is a principal solution at $x=0$ and $y_{-}$ is a non-principal solution at $x=0$.\\ 

\textit{Case i)} $\sqrt{-\mu}<\kappa_\al$. Notice that $y_{\pm}\in L^2_{1-\al}(0,1)$, thus $x=0$ is limit circle (LC), see also \cite[Theorem 7.2.2]{Zettl}. The result follows from Theorem 10.5.3 and equation (10.5.2) in \cite{Zettl} and by using Theorem 4.3 and equation (4.15) in \cite{Niessen}.\\

By using the notation in \cite[Chapter 10]{Zettl} we can see that $S_F=S_{\min}^*|_{D(S_F)}=S_{\max}|_{D(S_F)}$, and $S_{\max}u=w^{-1}Mu=\A u$, $u\in D_{\max}$.\\

\textit{Case ii)} $\sqrt{-\mu}\geq\kappa_\al$. Since $y_{-}\notin L^2_{1-\al}(0,1)$, then $x=0$ is limit point (LP). The result follows by using Theorem 10.4.4 in \cite{Zettl} with $A_1=1, A_2=0$.\\

The second part follows by using the computations in the proof of Proposition \ref{bases}. 
\end{proof}
\begin{remark}
From Theorem 10.5.3 and (10.5.2) in \cite{Zettl} we have that
\[\lim_{x\rightarrow 0^+}[u,y_+](x)=0\iff\lim_{x\rightarrow 0^+}\frac{u(x)}{y_-(x)}=0. \]
Notice that $py'_+y_-$ is constant on $(0,1)$. Therefore, the last condition is equivalent to 
\[\lim_{x\rightarrow 0^+}u'(x)py_+(x)=\lim_{x\rightarrow 0^+}x^{1+\sqrt{-\mu}}u'(x)=0,\quad u\in D(\A).\]
\end{remark}

Thus, we are in the same position as in \cite{GaloLopez}, so we can follow the same steps to get the proof of Theorem \ref{unosuma}.

\appendix
\section{Bessel functions}
We introduce the Bessel function of the first kind $J_{\nu}$ as follows
\begin{equation}\label{bessel}
J_{\nu}(x)=\sum_{m \geq 0} \frac{(-1)^{m}}{m ! \Gamma(m+\nu+1)}\left(\frac{x}{2}\right)^{2 m+\nu}, \quad x \geq 0,
\end{equation}
where $\Gamma(\cdot)$ is the Gamma function. In particular, for $\nu>-1$ and $0<x \leq \sqrt{\nu+1}$, from (\ref{bessel}) we have (see \cite[9.1.7, p. 360]{abram})
\begin{equation}\label{asincero}
J_{\nu}(x) \sim \frac{1}{\Gamma(\nu+1)}\left(\frac{x}{2}\right)^{\nu} \quad \text { as } \quad x \rightarrow 0^{+} .
\end{equation}
A Bessel function $J_\nu$ of the first kind solves the differential equation
\begin{equation}\label{Besselode}
x^2y''+xy'+(x^2-\nu^2)y=0.
\end{equation}
Bessel functions of the first kind satisfy the recurrence formula $([1], 9.1 .27)$:
\begin{equation}\label{recur}
x J_{\nu}^{\prime}(x)-\nu J_{\nu}(x)=-x J_{\nu+1}(x).
\end{equation}
Recall the asymptotic behavior of the Bessel function $J_{\nu}$ for large $x$, see \cite[Lem. 7.2, p. 129]{komo}.
\begin{lem}\label{asimxinf}
For any $\nu \in \mathbb{R}$
$$
J_{\nu}(x)=\sqrt{\frac{2}{\pi x}}\left\{\cos \left(x-\frac{\nu \pi}{2}-\frac{\pi}{4}\right)+\mathcal{O}\left(\frac{1}{x}\right)\right\} \quad \text { as } \quad x \rightarrow \infty
$$
\end{lem}

For $\nu >-1$ the Bessel function $J_{\nu}$ has an infinite number of real zeros $0<j_{\nu, 1}<j_{\nu, 2}<\ldots$, all of which are simple, with the possible exception of $x=0$. In \cite[Proposition 7.8]{komo} we can find the next information about the location of the zeros of the Bessel functions $J_{\nu}$:
\begin{lem}\label{consec}Let $\nu \geq 0$.\\
1. The difference sequence $\left(j_{\nu, k+1}-j_{\nu, k}\right)_{k}$ converges to $\pi$ as $k \rightarrow\infty$.\\
2. The sequence $\left(j_{\nu, k+1}-j_{\nu, k}\right)_{k}$ is strictly decreasing if $|\nu|>\frac{1}{2}$, strictly increasing if $|\nu|<\frac{1}{2}$, and constant if $|\nu|=\frac{1}{2}$.\\
\end{lem}

For $\nu \geq 0$ fixed, we consider the next asymptotic expansion of the zeros of the Bessel function $J_{\nu}$, see\cite[Section 15.53]{Watson},
\begin{equation}\label{asint}
j_{\nu, k}=\left(k+\frac{\nu}{2}-\frac{1}{4}\right) \pi-\frac{4 \nu^{2}-1}{8\left(k+\frac{\nu}{2}-\frac{1}{4}\right) \pi}+O\left(\frac{1}{k^{3}}\right), \quad \text { as } k \rightarrow\infty
\end{equation}

In particular, we have
\begin{equation}\label{below}
\begin{aligned}
&j_{\nu, k} \geq\left(k-\frac{1}{4}\right) \pi \quad \text { for } \nu \in\left[0, 1/2\right], \\
&j_{\nu, k} \geq\left(k-\frac{1}{8}\right) \pi \quad \text { for } \nu \in\left[1/2,\infty\right].
\end{aligned}
\end{equation}

\begin{lem}\label{reduce} For any $\nu \geq 0$ and any $k\geq 1$ we have
$$
\sqrt{j_{\nu, k}}\left|J_{\nu}^{\prime}\left(j_{\nu, k}\right)\right|=\sqrt{\frac{2}{\pi}}+O\left(\frac{1}{j_{\nu, k}}\right)\quad \text{as}\quad k \rightarrow \infty.
$$
\end{lem}
The proof of this result follows by using  (\ref{asincero}) and the recurrence formula (\ref{recur}).
\begin{lem} Let $\gamma=\gamma(\alpha,\beta,\mu)$ and $\nu=\nu(\alpha,\beta,\mu)$ given in (\ref{gamadef}) and (\ref{Nu}) respectively, then the $\al+\beta+\gamma$-generalized limit of $\Fi_k$ at $x=0$ is finite for all $k\geq 1$, and
	\begin{equation}\label{lim1}
	\OO_{\al+\beta+\gamma}(\Fi_k)= \frac{(2\kappa_\al)^{1/2}\pts{j_{\nu,k}}^{\nu}}{2^{\nu}\Gamma(\nu+1)\abs{J'_{\nu}\pts{j_{\nu, k}}}}, \quad k\geq1.
	\end{equation}
\end{lem}
\begin{proof}
	This result follows from (\ref{asincero}).
\end{proof}
We recall the following representation theorem, see \cite[p. 56]{koosis}.
\begin{thm}\label{repreTheo} Let $g(z)$ be an entire function of exponential type and assume that
	$$
	\int_{-\infty}^{\infty} \frac{\log ^{+}|g(x)|}{1+x^{2}} \mathrm{d}x<\infty.
	$$
	Let $\left\{b_{\ell}\right\}_{\ell \geq 1}$ be the set of zeros of $g(z)$ in the upper half plane $\Im(z)>0$ (each zero being repeated as many times as its multiplicity). Then,
	$$
	\log |g(z)|=A \Im(z)+\sum_{\ell=1}^{\infty} \log \left|\frac{z-b_{\ell}}{z-\bar{b}_{\ell}}\right|+\frac{\Im(z)}{\pi} \int_{-\infty}^{\infty} \frac{\log |g(s)|}{|s-z|^{2}} \mathrm{d}s,\quad\Im(z)>0,
	$$
	where
	$$
	A=\limsup _{y \rightarrow\infty} \frac{\log |g(i y)|}{y} .
	$$
\end{thm}


\begin{thebibliography}{99}
\bibitem{abram} Abramowitz M. and Stegun I. A., Handbook of mathematical functions with formulas, graphs and mathematical tables, National Bureau of Standards. App. Math. series, Vol. 55. 1964.
\bibitem{cazenave}Cazenave, Thierry; Haraux, Alain, An introduction to semilinear evolution equations.
Oxford Lecture Series in Mathematics and its Applications, 13. The Clarendon Press, Oxford University Press, New York, 1998. xiv+186 pp. ISBN: 0-19-850277-X
\bibitem{dautray}Dautray, Robert; Lions, Jacques-Louis, Mathematical analysis and numerical methods for science and technology. Vol. 2. Functional and variational methods. Springer-Verlag, Berlin, 1988. xvi+561 pp. ISBN: 3-540-19045-7.
\bibitem{du}Du, Runmei Null controllability for a class of degenerate parabolic equations with the gradient terms. J. Evol. Equ. 19 (2019), no. 2, 585--613.
\bibitem{flores}Flores, C. and Teresa, L.; Carleman estimates for degenerate parabolic equations with first order
terms and applications, C. R. Math. Acad. Sci. Paris, 348(2010), 391--396.
\bibitem{flores2} Flores, C. and Teresa, L.; Null controllability of one-dimensional degenerate parabolic equations with first-order terms. Discrete and Continuous Dynamical Systems - B, 2020, 25 (10): 3963--3981.
\bibitem{gueye2}Gueye, Mamadou; Lissy, P.; Singular optimal control of a 1-D parabolic-hyperbolic degenerate equation. ESAIM Control Optim. Calc. Var. 22 (2016), no. 4, 1184--1203.
\bibitem{Hoch}Hochstadt, Harry; The mean convergence of Fourier-Bessel series. SIAM Rev. 9 (1967), 211--218.
\bibitem{koosis}Koosis P., The logarithmic integral I \& II, Cambridge Studies in Advanced Mathematics 12 (1988) \& Cambridge Studies in Advanced Mathematics 21 (1992), Cambridge University Press, Cambridge.
\bibitem{komo}Komornik V. and Loreti P., Fourier series in control theory. Springer (2005).
\bibitem{Fatorini}Fattorini H.O. and Russell D.L., Exact controllability theorems for linear parabolic equations in one space dimension. Arch. Ration. Mech. Anal. 43 (1971) 272--292. 
\bibitem{GaloLopez}L. Galo and M. López-García, Boundary controllability for a 1D degenerate parabolic equation with drift and a singular potential. arXiv:2302.01197
\bibitem{Niessen} Niessen H., Zettl A., Singular Sturm-Liouville problems: The Friedrichs extension and Comparison of eigenvalues, Proc. London Math. Soc. v.64 (1992) 545-578.
\bibitem{Tucsnak}Tenenbaum G. and Tucsnak M.; New blow-up rates for fast controls of Schrodinger and heat equations. J. Differ. Equ. 243 (2007) 70–100.
\bibitem{Watson}Watson G.N., A treatise on the theory of Bessel functions. Cambridge University Press, Cambridge, England (1958). 
\bibitem{Zettl} Zettl A., Sturm-Liouville Theory, Mathematical Surveys and Monographs, vol. 121. Am. Math. Soc., Providence (2005).
\end{thebibliography}
\end{document}